\documentclass[11pt,a4paper]{amsart}
\usepackage{cite}
\usepackage[T1]{fontenc}
\usepackage{lmodern}
\usepackage{amsmath}
\usepackage{amsfonts}
\usepackage{amssymb}
\usepackage{mathtools}
\usepackage{latexsym}
\usepackage{parskip}
\usepackage{enumerate}
\usepackage{lipsum}
\usepackage[latin1]{inputenc}
\usepackage{amsmath}
\usepackage{amsfonts}
\usepackage{amssymb}
\usepackage{latexsym}

\usepackage{verbatim}

\usepackage{multicol}

\usepackage{parskip}
\usepackage{fullpage}
\usepackage[pdftex]{graphicx}

\usepackage{soul}

\usepackage{enumerate}
\numberwithin{equation}{section}

\usepackage{amsthm}

\DeclareMathOperator{\e}{e}

\newtheorem{theorem}{Theorem}[section]
\newtheorem{corollary}[theorem]{Corollary}

\theoremstyle{definition}
\newtheorem{definition}[theorem]{Definition}
\newtheorem{remark}[theorem]{Remark}
\newtheorem{example}[theorem]{Example}

\newcommand{\Id}{\mathrm{Id}}

\usepackage[usenames,dvipsnames,svgnames,table]{xcolor}

\DeclareMathOperator{\Ker}{Ker}

\newcommand{\N}{\mathbb{N}}
\newcommand{\R}{\mathbb{R}}

\title{Generalized dichotomies via time rescaling}
\author{Davor Dragi\v{c}evi\'c}
\address{Faculty of Mathematics, University of Rijeka, Croatia}
\email[Davor Dragi\v cevi\' c]{ddragicevic@math.uniri.hr}
\author{C\'{e}sar M. Silva}
\address{Centro de Matem\'atica e Aplica\c{c}\~oes\\
Universidade da Beira Interior\\
6201-001 Covilh\~a\\
Portugal}
\email[C\'{e}sar M. Silva]{csilva@ubi.pt}
\begin{document}

\maketitle

\begin{abstract}
For discrete-time nonautonomous linear dynamics and a large class of discrete growth rates $\mu$, we show that the notion of $\mu$ dichotomy (with respect to a sequence of norms) can be completely characterized in terms of ordinary and exponential dichotomy (with respect to a sequence of norms) by employing a suitable rescaling of time. Previously, such a result was known only in the particular case of polynomial dichotomies.
As a nontrivial application of our results, we study the structure of a generalized Sacker-Sell spectrum and obtain certain nonautonomous topological and smooth linearization results.
\end{abstract}


\section{Introduction}

The notion of exponential dichotomy introduced by Perron~\cite{Per} plays an important role in the qualitative study of nonautonomous dynamical systems. 	It corresponds to assuming the existence of an exponential contraction and
expansion along complementary directions at each moment of time, thus representing a natural counterpart to the notion of hyperbolicity for autonomous dynamics.  Among many important consequences, we mention the existence of stable and unstable manifolds as well as topological linearization for small nonlinear perturbations of linear dynamics exhibiting an exponential dichotomy. We refer to~\cite{BDV, CL, Coppel, DK, Hen, KR, MS2, Pot, SY} for a detailed discussion of various aspects of this theory.

Despite its importance, in some situations, the exponential dichotomy might be regarded
restrictive since it requires that the rates of contraction and expansion
along the stable and unstable directions be exponential. Indeed, it is fairly easy to construct examples of nonautonomous dynamics which admit a splitting
into stable and unstable directions but with non-exponential rates of contraction and expansion along these directions. To our knowledge, Muldowney~\cite{M} and Naulin and Pinto~\cite{NP}
were the first to study non-exponential dichotomies. More recently, Barreira and Valls~\cite{BV} initiated a systematic study of such dichotomies in the nonuniform framework connecting the existence of this type of behavior with the nonvanishing of certain generalized Lyapunov exponents (see also~\cite{BS2, BS1} for a subsequent work by Bento and Silva).   We stress that particular attention has been paid to the so-called polynomial dichotomies~\cite{BV-09-01, BS}. We refer to~\cite{BD,BLMS,BMT,Chu,DLP,JZF,LP1,LP2} and references therein for various results devoted to dichotomies with growth rates.

Building on earlier work on exponential~\cite{BDV1} and polynomial dichotomies~\cite{D1}, Silva~\cite{Silva} introduced, for a nonautonomous linear dynamics with discrete time, the notion of $\mu$-dichotomy with respect to a sequence of norms, where $\mu$ is a discrete growth rate. We emphasize that this notion includes the notion of a nonuniform $\mu$-dichotomy as a particular case (see~\cite[Theorem 4.1]{Silva}).

In order to describe the results of the present paper, let us consider a nonautonomous linear difference equation
\begin{equation}\label{LDE-intro}
    x_{n+1}=A_nx_n \quad n\in \N,
\end{equation}
on an arbitrary Banach space $X$, where $\{A_n \}_{n\in \N}$ is a sequence of bounded linear operators on $X$. It follows from our first main result (see Theorem~\ref{teo:main}) that (provided that $\mu$ is slowly varying) the notion of $\mu$-dichotomy with respect to a sequence of norms for~\eqref{LDE-intro} can be completely characterized in terms of the notion of ordinary and exponential dichotomy with respect to a sequence of norms for~\eqref{LDE-intro} and a system obtained from~\eqref{LDE-intro} by a suitable rescaling of time, respectively. Moreover, we show (see Theorem~\ref{newthm}) that under certain additional conditions such a characterization holds without involving the existence of an ordinary dichotomy for~\eqref{LDE-intro}. These results show that, for a large class of growth rates $\mu$,  $\mu$-dichotomies with respect to a sequence of norms can be described as exponential dichotomies with respect to a sequence of norms after a suitable time rescaling.  We stress that this was previously known only in the particular case of polynomial dichotomies~\cite{DSS1} (see also~\cite{BDZ,DSS2} for some applications).

We proceed by giving several nontrivial applications of these results. Firstly, we show that in some cases the version of the Sacker-Sell spectrum for~\eqref{LDE-intro} introduced with respect to the notion of (uniform) $\mu$-dichotomy coincides with the classical Sacker-Sell spectrum of a system obtained from~\eqref{LDE-intro} via suitable time rescaling (see Theorem~\ref{dicspec}). 

Next, we consider nonlinear perturbations of~\eqref{LDE-intro} of the form
\begin{equation}\label{NDE-intro}
    x_{n+1}=A_n x_n +f_n(x_n) \quad n\in \N,
\end{equation}
where $\{f_n\}_{n\in \N}$ is a sequence of maps $f_n\colon X\to X$. We obtain a series of results devoted to the linearization of~\eqref{NDE-intro}. 
More precisely, we show that if~\eqref{LDE-intro} admits a strong nonuniform $\mu$-dichotomy and the nonlinearities $f_n$ are ``small'' in the appropriate sense, then~\eqref{LDE-intro} and~\eqref{NDE-intro} are topologically equivalent and the conjugacies are locally H\"{o}lder continuous (see Theorem~\ref{NTL}). We stress that a similar result was established (using different techniques) by Barreira and Valls~\cite{BV-09}. However, as explained in Remark~\ref{comparison} there are important differences between the two and our Theorem~\ref{NTL} seems to be the first result which is applicable to the case of polynomial dichotomies (either uniform or nonuniform). We note that the research devoted to topological nonautonomous linearization was initiated by Palmer~\cite{Palmer}.

Next, we obtain a Sternberg-type linearization result (see Theorem~\ref{Sieg}) that formulates a $C^\ell$-linearization result for~\eqref{NDE-intro}. It works under the assumptions that~\eqref{LDE-intro} exhibits a (uniform) $\mu$-dichotomy, that there are no resonances up to a certain order (formulated in terms of the generalized Sacker-Sell spectrum) and that nonlinearities $f_n$ are ``small'' in an appropriate sense.  We deduce Theorem~\ref{Sieg} as a consequence of the main results described above of our paper and the Sternberg-type theorem obtained in~\cite{CDS} which considers~\eqref{LDE-intro} admitting the (uniform) exponential dichotomy. For other results dealing with nonautonomous smooth linearization we refer to~\cite{DZZ1, DZZ2}.

\section{Time rescaling}
	Throughout this paper $X=(X, \| \cdot \|)$ is an arbitrary Banach space. By $\mathcal B(X)$ we denote the space of all bounded linear operators on $X$ equipped with the operator norm, which we also denote by $\| \cdot \|$.
	Let $\{A_n\}_{n \in \mathbb{N}} \subset \mathcal{B}(X)$ and consider the nonautonomous linear system
	\begin{equation}\label{eq:Lin}
	x_{n+1}=A_n x_n, \ \ n \in \mathbb{N}. \tag{A}
	\end{equation}
	Let $\Phi_A=\left\{\Phi_A(m, k)\right\}_{m \ge k}$ be the evolution family associated to~\eqref{eq:Lin} defined by
	\begin{equation}\label{eq:evolution-family-Lin}
	\Phi_A(m, k)
	=\left\{\begin{array}{ll}
		A_{m-1} \cdots A_k, & m>k \\
		\Id, & m=k.
	\end{array} \right.
	\end{equation}
	Here $\Id$ denotes the identity operator on $X$.
We say that a sequence $\mu=(\mu_n)_{n \in \N_0}$ is a discrete growth rate if it is positive, strictly increasing, and satisfies $\lim\limits_{n\to \infty} \mu_n = +\infty$ and $\mu_0=1$. 

Given a discrete growth rate $\mu=(\mu_n)_{n \in \mathbb{N}_0}$, we can associate to it the continuous, strictly increasing function $\widetilde \mu: [0,+\infty) \to [1,+\infty )$ given by
\[
\widetilde \mu(t)=
\begin{cases}
	\mu_n & \text{ if } \ t=n\\
	\mu_n+(t-n)(\mu_{n+1}-\mu_n) & \text{ if } \ n<t<n+1
\end{cases}, \quad n \in \N_0.
\]
Clearly, this is an invertible function and its inverse ${\widetilde \mu}^{-1}:[1,+\infty) \to [0,+\infty)$ is the continuous, strictly increasing function given by
\[
\widetilde \mu^{-1}(t)=
\begin{cases}
	n & \text{ if } \ t=\mu_n\\
	n+\dfrac{t-\mu_n}{\mu_{n+1}-\mu_n} & \text{ if } \ \mu_n<t<\mu_{n+1}
\end{cases}, \quad n \in \N_0.
\]

\begin{definition}\label{def:ordinary-dich}
	We say that~\eqref{eq:Lin} (or the sequence $\{A_n\}_{n\in \N}$) admits an ordinary dichotomy with respect to a sequence of norms $\left\{\|\cdot\|_k\right\}_{k \in \mathbb{N}}$ on $X$
 if there are a sequence of projections $\{P_k\}_{k \in \mathbb{N}}$ on $X$ and $K \geq 1$ such that:
	\begin{enumerate}[(od1)]
		\item\label{def:od1} $A_k P_k=P_{k+1} A_k$, for all $k \in \mathbb{N}$;
		\item\label{def:od3} for every $k \in \mathbb{N}, {A_k}_{\mid}: \operatorname{Ker} P_k \rightarrow \operatorname{Ker} P_{k+1}$ is an isomorphism, where ${A_k}_{\mid}$ denotes the restriction of $A_k$ to $\Ker P_k$;
		\item\label{def:od4} $\left\|\Phi_A(m, k) P_k x\right\|_m \leq K\|x\|_k$, for all $x \in X$ and all $m \ge k$;
		\item\label{def:od5} $\left\|\Phi_A(m, k) (\Id-P_k)x\right\|_m \leq K\|x\|_k$, for all $x \in X$ and all $m\le k$.
	\end{enumerate}
\end{definition}

We remark that~(od\ref{def:od1}) is equivalent to
\begin{equation}\label{eq:compatibility-evol-family}
	\Phi_A(m, n) P_n=P_m \Phi_A(m,n), \quad \text{for all $m \ge  n \ge 1$}
\end{equation}
and~(od\ref{def:od3}) is equivalent to the invertibility of $\Phi_A(m,n)_{\mid}: \operatorname{Ker} P_n \rightarrow \operatorname{Ker} P_m$, for all $m \ge n \ge 1$.
Note that, in~(od\ref{def:od5}), $\Phi_A(m,k)$ denotes the inverse of the restriction of $\Phi_A(k,m)$ to $\Ker P_m$.

\begin{definition}\label{def:mu-dich-norms}
Let $\mu=(\mu_n)_{n\in \N_0}$ be a discrete growth rate.
We say that~\eqref{eq:Lin} (or the sequence $\{A_n\}_{n\in \N}$) admits a $\mu$-dichotomy with respect to a sequence of norms $\left\{\|\cdot\|_k\right\}_{k \in \mathbb{N}}$ on $X$
  if there are $N \geq 1$, $\nu>0$ and a sequence of projections $\{P_k\}_{k \in \mathbb{N}}$ on $X$ such that (od\ref{def:od1})--(od\ref{def:od3}) in Definition~\ref{def:ordinary-dich} are satisfied and additionally:
\begin{enumerate}[($\mu$1)]
	\item\label{def:mu1} 
	$\left\|\Phi_A(m, k) P_k x\right\|_m \leq N \left(\frac{\mu_m}{\mu_k}\right)^{-\nu}\|x\|_k$, for every $x \in X$ and $m \ge k$;
	\item\label{def:mu2} $\left\|\Phi_A(m, k) (\Id-P_k)x\right\|_m \leq N \left(\frac{\mu_k}{\mu_m}\right)^{-\nu}\|x\|_k$, for every $x \in X$ and  $m \le k$.
\end{enumerate}
\end{definition}

\begin{remark}
When, in Definition~\ref{def:mu-dich-norms}, we set $\mu_n=\e^n$, for $n \in \N_0$, we say that we have an exponential dichotomy with respect to a sequence of norms $\left\{\|\cdot\|_k\right\}_{k \in \mathbb{N}}$. Similarly, when we let $\mu_n=n+1$ for $n \in \N_0$, we say that we have a polynomial dichotomy with respect to a sequence of norms $\left\{\|\cdot\|_k\right\}_{k \in \mathbb{N}}$.
\end{remark}

Given discrete growth rates $\mu=(\mu_n)_{n \in\N_0}$ and $\eta=(\eta_n)_{n \in\N_0}$, we associate with the system~\eqref{eq:Lin} the family of linear operators $Q^{\mu,\eta}=\{Q^{\mu,\eta}_n\}_{n \in \mathbb{N}} \subset \mathcal{B}(X)$, given by 
\begin{equation}\label{eq:Q-mu-eta}
Q^{\mu,\eta}_n=\Phi_A\left(\lfloor \widetilde\mu^{-1}(\eta_n) \rfloor +1, \  \lfloor \widetilde\mu^{-1}(\eta_{n-1}) \rfloor +1 \right), \quad \forall n \in \mathbb{N}.
\end{equation}
Consider the nonautonomous linear system
\begin{equation}\label{eq:Lin-Q}
	y_{n+1}=Q^{\mu,\eta}_n \, y_n, \quad n \in \mathbb{N}. \tag{$Q^{\mu,\eta}$}
\end{equation}
It is easy to obtain the evolution family $\Phi_{Q^{\mu,\eta}}=\left\{\Phi_{Q^{\mu,\eta}}(m,n)\right\}_{m \ge n}$ associated to~\eqref{eq:Lin-Q} in terms of the evolution family $\Phi_A=\left\{\Phi_A(m,n)\right\}_{m \ge n}$ associated to~\eqref{eq:Lin}:
\[
\begin{split}
\Phi_{Q^{\mu,\eta}}(m,n)
& =\prod_{k=n}^{m-1}Q^{\mu,\eta}_k=\prod_{k=n}^{m-1}\Phi_A(\lfloor\widetilde\mu^{-1}(\eta_{k})\rfloor +1,\,\lfloor\widetilde\mu^{-1}(\eta_{k-1})\rfloor +1)\\
& =\Phi_A(\lfloor\widetilde\mu^{-1}(\eta_{m-1})\rfloor +1,\,\lfloor\widetilde\mu^{-1}(\eta_{n-1})\rfloor +1), \quad m\ge n.
\end{split}
\]

Given a sequence of norms $\left\{\|\cdot\|_k\right\}_{k \in \mathbb{N}}$ and a growth rate $\eta$, we consider a new sequence of norms $\left\{\|\cdot\|_k^\eta\right\}_{k \in \mathbb{N}}$ given by
$$
\|\cdot\|_k^\eta:=\|\cdot\|_{\lfloor\widetilde\mu^{-1}(\eta_{k-1})\rfloor+1}, \quad k\in \N.
$$
Moreover, given a linear operator $T \colon X\to X$ we let, for all $k \in \N$, 
\[
\begin{split}
\|T\|_k^\eta=
\inf \{M \in (0,+\infty]: \|Tx\|_{\lfloor\widetilde\mu^{-1}(\eta_{k-1})\rfloor+1}\le M\|x\|_{\lfloor\widetilde\mu^{-1}(\eta_{k-1})\rfloor+1}, \ \forall x \in X \}.
\end{split}
\] We note that, provided that $\| \cdot \|_{\lfloor\widetilde\mu^{-1}(\eta_{k-1})\rfloor+1}$ is equivalent to $\|\cdot \|$, we have $\|T\|_k^\eta<+\infty$ for every $T\in \mathcal B(X)$.
The central result of this section is the following.

\begin{theorem}\label{teo:main}
Let $\left \{\| \cdot \|_k \right \}_{k\in \N}$ be a sequence of norms on $X$ and  $\mu$ a growth rate satisfying 
\begin{equation}\label{eq:bound-rates-1-a}
	\frac{\mu_{n+1}}{\mu_n}\le \theta, \quad \text{for all} \quad n \in \N_0,
\end{equation} 
for some $\theta\ge 1$. The following assertions are equivalent:
\begin{enumerate}[(i)]
\item\label{teo:main-1} system~\eqref{eq:Lin} admits a $\mu$-dichotomy with respect to the sequence of norms $\left\{\|\cdot\|_k\right\}_{k \in \mathbb{N}}$ relative to a sequence of projections $\{P_k\}_{k \in \N}$;
\item\label{teo:main-2} system~\eqref{eq:Lin} admits an ordinary dichotomy with respect to the sequence of norms $\left\{\|\cdot\|_k\right\}_{k \in \mathbb{N}}$ relative to a sequence of projections $\{P_k\}_{k \in \mathbb{N}}$ and, for any growth rate $\eta$ satisfying 
\begin{equation}\label{eq:bound-rates-2}
	\frac{\eta_{n+1}}{\eta_n} \le \theta, \quad \text{for all} \quad n \in \N_0,
\end{equation}
for some $\theta\ge 1$, system~\eqref{eq:Lin-Q} admits a $\eta$-dichotomy with respect to the norms $\left\{\|\cdot\|_k^{\eta}\right\}_{k \in \mathbb{N}}$ relative to the sequence of projections $\left\{P^\eta_k\right\}_{k \in \mathbb{N}}$, where $P^\eta_k=P_{\lfloor \widetilde{\mu}^{-1}(\eta_{k-1})\rfloor +1}$, for all $k \in \mathbb{N}$;
\item\label{teo:main-3} system~\eqref{eq:Lin} admits an ordinary dichotomy with respect to the sequence of norms $\left\{\|\cdot\|_k\right\}_{k \in \mathbb{N}}$ relative to a sequence of projections $\{P_k\}_{k \in \mathbb{N}}$ and there is a growth rate $\eta$ satisfying~\eqref{eq:bound-rates-2} such that system~\eqref{eq:Lin-Q} admits a $\eta$-dichotomy with respect to the sequence of norms $\left\{\|\cdot\|_k^{\eta}\right\}_{k \in \mathbb{N}}$ relative to the sequence of projections $\left\{P^\eta_k\right\}_{k \in \mathbb{N}}$, where $P^\eta_k=P_{\lfloor \widetilde{\mu}^{-1}(\eta_{k-1})\rfloor +1}$, for all $k \in \mathbb{N}$.
\end{enumerate} 
\end{theorem}

\begin{proof}
$\eqref{teo:main-1} \Longrightarrow \eqref{teo:main-2} \quad$ Assume that the system~\eqref{eq:Lin} admits a $\mu$-dichotomy with respect to the sequence of norms $\left\{\|\cdot\|_k\right\}_{k \in \mathbb{N}}$ relative to a sequence of projections $\{P_k\}_{k \in \mathbb{N}}$ and let $N \geq 1, \nu>0$ be given by Definition~\ref{def:mu-dich-norms}. Then, for all $x \in X$ and all $m \ge k$, we have
\begin{equation}\label{mu1}\left\|\Phi_A(m, k) P_k x\right\|_m \leq N \left(\frac{\mu_m}{\mu_k}\right)^{-\nu}\|x\|_k\end{equation}
and, for all $x \in X$ and all $m \le k$, we have
\begin{equation}\label{mu2}\left\|\Phi_A(m, k) (\Id-P_k)x\right\|_m \leq N \left(\frac{\mu_k}{\mu_m}\right)^{-\nu}\|x\|_k.\end{equation}

Since $\frac{\mu_r}{\mu_s}\ge 1$ for $r \ge s$, we immediately conclude that the system~\eqref{eq:Lin} admits an ordinary dichotomy with respect to the sequence of norms $\left\{\|\cdot\|_k\right\}_{k \in \mathbb{N}}$ and the same sequence of projections $\{P_k\}_{k \in \mathbb{N}}$.


On the other hand, by~\eqref{eq:compatibility-evol-family}, we obtain, for all $k \in \mathbb{N}$,
\[
\begin{split}
Q^{\mu,\eta}_k P^\eta_k 
& =\Phi_A\left(\lfloor \widetilde\mu^{-1}(\eta_k) \rfloor +1, \  \lfloor \widetilde\mu^{-1}(\eta_{k-1}) \rfloor +1\right) P_{\lfloor \widetilde{\mu}^{-1}(\eta_{k-1})\rfloor +1} \\
& =P_{\lfloor \widetilde{\mu}^{-1}(\eta_{k})\rfloor +1} \Phi_A\left(\lfloor \widetilde\mu^{-1}(\eta_k) \rfloor +1, \  \lfloor \widetilde\mu^{-1}(\eta_{k-1}) \rfloor +1\right) =P^\eta_{k+1} Q^{\mu,\eta}_k.
\end{split}
\]

By~(od\ref{def:od3}) we conclude that, for every $k \in \mathbb{N}$, $$\Phi_A\left(\lfloor \widetilde\mu^{-1}(\eta_k) \rfloor +1, \  \lfloor \widetilde\mu^{-1}(\eta_{k-1}) \rfloor +1\right)_\mid: \operatorname{Ker} P_{\lfloor \widetilde\mu^{-1}(\eta_{k-1}) \rfloor +1} \rightarrow \operatorname{Ker} P_{\lfloor \widetilde\mu^{-1}(\eta_k) \rfloor +1}$$ 
is invertible. In other words, ${Q^{\mu,\eta}_k}_{\mid}: \operatorname{Ker} P^\eta_k \rightarrow \operatorname{Ker} P^\eta_{k+1}$ is invertible.

Note that, by~\eqref{eq:bound-rates-1-a}, for each $t \in \N_0$, we have
\[
\frac{\tilde{\mu}(t+1)}{\tilde{\mu}(t)}=\frac{\mu_{t+1}}{\mu_t}\le \theta
\]
and, for each $t \in [0,+\infty)\setminus\N_0$, we have
\[
\frac{\tilde{\mu}(t+1)}{\tilde{\mu}(t)}=\frac{\mu_{r+1}+(t-r)(\mu_{r+2}-\mu_{r+1})}{\mu_r+(t-r)(\mu_{r+1}-\mu_r)}\le \frac{\mu_{r+2}}{\mu_r}=\frac{\mu_{r+2}}{\mu_{r+1}}\frac{\mu_{r+1}}{\mu_r} \le \theta^2,
\]
where $r=\lfloor t\rfloor \in \N_0$. We conclude that
\begin{equation}\label{eq:bound2}
\frac{\tilde{\mu}(t+1)}{\tilde{\mu}(t)}\le \theta^2, \quad \text{for all} \ t \in [0,+\infty).
\end{equation}
Let $k \in \mathbb{N}$ and $x \in X$. Using \eqref{eq:bound-rates-1-a}, \eqref{eq:bound-rates-2}, \eqref{mu1} and~\eqref{eq:bound2}, we have, for all $m \ge k \ge 1$, 
\[
\begin{split}
	& \quad \left\|\Phi_{Q^{\mu,\eta}}(m,k)P_k^\eta x\right\|_m^{\eta}\\
    & =\left\|\Phi_A(\lfloor\widetilde\mu^{-1}(\eta_{m-1})\rfloor +1,\,\lfloor\widetilde\mu^{-1}(\eta_{k-1})\rfloor+1)P_{\lfloor\widetilde\mu^{-1}(\eta_{k-1})\rfloor +1 }x \right\|_{\lfloor\widetilde\mu^{-1}(\eta_{m-1})\rfloor +1}\\
  & \leq N \left(\frac{\mu_{\lfloor\widetilde\mu^{-1}(\eta_{m-1})\rfloor +1}}{\mu_{\lfloor\widetilde\mu^{-1}(\eta_{k-1})\rfloor +1}}\right)^{-\nu}\|x\|_{\lfloor\widetilde\mu^{-1}(\eta_{k-1})\rfloor +1}\\
  &\leq N \theta^\nu \left(\frac{\mu_{\lfloor\widetilde\mu^{-1}(\eta_{m-1})\rfloor }}{\mu_{\lfloor\widetilde\mu^{-1}(\eta_{k-1})\rfloor }}\right)^{-\nu}\|x\|_{\lfloor\widetilde\mu^{-1}(\eta_{k-1})\rfloor +1}\\
  & = N \theta^\nu \left(\frac{\widetilde\mu(\lfloor\widetilde\mu^{-1}(\eta_{m-1})\rfloor)}{\widetilde\mu(\widetilde\mu^{-1}(\eta_{m-1}))} \
	\frac{\widetilde\mu(\widetilde\mu^{-1}(\eta_{m-1}))}{\widetilde\mu(\widetilde\mu^{-1}(\eta_{m}))} \ \frac{\widetilde\mu(\widetilde\mu^{-1}(\eta_{m}))}{\widetilde\mu(\widetilde\mu^{-1}(\eta_{k}))} \ \frac{\widetilde\mu(\widetilde\mu^{-1}(\eta_{k}))}{\widetilde\mu(\widetilde\mu^{-1}(\eta_{k-1}))} \ \frac{\widetilde\mu(\widetilde\mu^{-1}(\eta_{k-1}))}{\widetilde\mu(\lfloor\widetilde\mu^{-1}(\eta_{k-1})\rfloor)}\right)^{-\nu} \|x\|_k^\eta \\
    & \leq N \theta^\nu \left(\frac{\widetilde\mu(\widetilde\mu^{-1}(\eta_{m-1}))}{\widetilde\mu(\widetilde\mu^{-1}(\eta_{m-1})-1)} \ \frac{\eta_{m}}{\eta_{m-1}}\ \frac{\eta_{k}}{\eta_{m}} \ \frac{\eta_{k-1}}{\eta_{k}} \ \frac{\widetilde\mu(\lfloor\widetilde\mu^{-1}(\eta_{k-1})\rfloor)}{\widetilde\mu(\widetilde\mu^{-1}(\eta_{k-1}))}\right)^\nu \|x\|_k^\eta\\
   & \leq N\theta^{4\nu} \left(\frac{\eta_{m}}{\eta_k}\right)^{-\nu}\|x\|_k^\eta,
\end{split}
\]
with the convention that  $\widetilde{\mu}^{-1}(\eta_{m-1})-1:=0$  if $\widetilde{\mu}^{-1}(\eta_{m-1})<1$.

Let $k \in \mathbb{N}$ and $x \in X$.  Using~ \eqref{eq:bound-rates-2}, \eqref{mu2} and~\eqref{eq:bound2}, we have, for all $1\le m \le k$,
$$
\begin{aligned}
	& \quad \left\|\Phi_{Q^{\mu,\eta}}(m,k)(\Id-P^\eta_k)x\right\|_m^{\eta}\\
	& =\left\|\Phi_A(\lfloor\widetilde\mu^{-1}(\eta_{m-1})\rfloor +1,\,\lfloor\widetilde\mu^{-1}(\eta_{k-1})\rfloor+1)(\Id-P_{\lfloor\widetilde\mu^{-1}(\eta_{k-1})\rfloor+1})x\right\|_{\lfloor\widetilde\mu^{-1}(\eta_{m-1})\rfloor+1}\\
	& \leq N \left(\frac{\mu_{\lfloor\widetilde\mu^{-1}(\eta_{k-1})\rfloor+1}}{\mu_{\lfloor\widetilde\mu^{-1}(\eta_{m-1})\rfloor+1}}\right)^{-\nu} \|x\|_{\lfloor\widetilde\mu^{-1}(\eta_{k-1})\rfloor+1}\\
&\leq N\theta^\nu \left(\frac{\mu_{\lfloor\widetilde\mu^{-1}(\eta_{k-1})\rfloor}}{\mu_{\lfloor\widetilde\mu^{-1}(\eta_{m-1})\rfloor}}\right)^{-\nu} \|x\|_{\lfloor\widetilde\mu^{-1}(\eta_{k-1})\rfloor+1}\\
	& \leq N\theta^\nu
	 \left(\frac{\eta_{m-1}}{\eta_{m}} \ \frac{\eta_m}{\eta_k} \ \frac{\eta_k}{\eta_{k-1}} \ \frac{\widetilde{\mu}(\widetilde{\mu}^{-1}(\eta_{k-1})}{\widetilde{\mu}(\widetilde{\mu}^{-1}(\eta_{k-1}-1))}\right)^{\nu} \|x\|_k^\eta\\
	 & \leq N\theta^{4\nu}
	 \left(\frac{\eta_k}{\eta_m} \right)^{-\nu} \|x\|_k^\eta.	 
\end{aligned}
$$

It follows that~\eqref{eq:Lin-Q} admits an $\eta$-dichotomy with respect to the sequence of norms $\left\{\|\cdot\|_k^{\eta}\right\}_{k \in \mathbb{N}}$ relative to the sequence of projections $\left\{P^\eta_k\right\}_{k \in \mathbb{N}}$.

$\eqref{teo:main-2} \Longrightarrow \eqref{teo:main-3} \ \ $ This implication is immediate.

$\eqref{teo:main-3} \Longrightarrow \eqref{teo:main-1} \quad$ Assume that~\eqref{eq:Lin} admits an ordinary dichotomy with respect to the sequence of norms $\left\{\|\cdot\|_k\right\}_{k \in \mathbb{N}}$ relative to a sequence of projections $\{P_k\}_{k \in \mathbb{N}}$: for some $K \geq 1$, we have, for all $x\in X$ and $m \ge k$,
\begin{equation}\label{eq:ordinary-dich1-iii->i}
\left\|\Phi_A(m, k) P_k x\right\|_m \leq K\|x\|_k
\end{equation}
and, for all $x\in X$ and $m \le k$,
\begin{equation}\label{eq:ordinary-dich2-iii->i}
\left\|\Phi_A(m, k) (\Id-P_k)x\right\|_m \leq K\|x\|_k. 
\end{equation}

Assume also that system~\eqref{eq:Lin-Q} admits a $\eta$-dichotomy with respect to the sequence of norms $\left\{\|\cdot\|_k^\eta\right\}_{k \in \mathbb{N}}$ relative to the projections $\left\{P^\eta_k\right\}_{k \in \mathbb{N}}$: there are $N \geq 1, \nu>0$ such that, for all $x \in X$ and $m \ge k$, we have
\begin{equation}\label{eq:eta-dich1-iii->i}
\left\|\Phi_{Q^{\mu,\eta}}(m,k) P_k^\eta x\right\|_m^\eta \leq N \left(\frac{\eta_m}{\eta_k}\right)^{-\nu}\|x\|_k^\eta
\end{equation}
and for $m\le k$,
\begin{equation}\label{eq:eta-dich2-iii->i}
\left\|\Phi_{Q^{\mu,\eta}}(m,k) (\Id-P_k^\eta)x\right\|_m^\eta \leq N \left(\frac{\eta_k}{\eta_m}\right)^{-\nu}\|x\|_k^\eta.
\end{equation}
We claim that there exists $N_1\ge 1$ such that  for $x \in X$ and $m \ge k \ge 1$, we have
\begin{equation}\label{eq:claim1}
	\left\|\Phi_A(m, k) P_k x\right\|_m
	\le N_1 \left(\frac{\mu_{m}}{\mu_{k}}\right)^{-\nu} \left\|x\right\|_k.
\end{equation}
We split the proof of~\eqref{eq:claim1} into two cases:\\
\underline{Case I} - There are $n,j \in \N_0$ with $n>j$ such that 
$$\lfloor\widetilde\mu^{-1}(\eta_{n+1})\rfloor +1> m \ge \lfloor\widetilde\mu^{-1}(\eta_n)\rfloor +1 \ge \lfloor\widetilde\mu^{-1}(\eta_{j+1})\rfloor +1 >  k
\ge \lfloor\widetilde\mu^{-1}(\eta_j)\rfloor+1.$$
 Recalling that $P_{\lfloor {\widetilde \mu}^{-1} (\eta_{j+1}) \rfloor +1}=P_{j+2}^\eta$ and $\|\cdot\|^\eta_{j+2}=\|\cdot\|_{\lfloor {\widetilde \mu}^{-1} (\eta_{j+1}) \rfloor+1}$, by~\eqref{eq:ordinary-dich1-iii->i} and~\eqref{eq:eta-dich1-iii->i} we get 
\[
\begin{split}
& \left\|\Phi_A(m, k) P_k x\right\|_m \\
	& =\left\|
	\Phi_A\left(m, \lfloor\widetilde\mu^{-1}(\eta_n)\rfloor +1\right) P_{\lfloor \widetilde{\mu}^{-1}(\eta_n)\rfloor+1} \Phi_A\left(\lfloor\widetilde\mu^{-1}(\eta_n)\rfloor +1, k\right) P_k x\right\|_m \\
	& \le K \left\| \Phi_A\left(\lfloor\widetilde\mu^{-1}(\eta_n)\rfloor +1, \lfloor\widetilde\mu^{-1}(\eta_{j+1})\rfloor +1\right)
	\Phi_A\left(\lfloor\widetilde\mu^{-1}(\eta_{j+1})\rfloor +1,k\right) P_k x\right\|_{\lfloor\widetilde\mu^{-1}(\eta_n)\rfloor +1}\\
&=K\left\|\Phi_{Q^{\mu,\eta}}\left(n+1,j+2\right)P_{j+2}^\eta \Phi_A\left(\lfloor\widetilde\mu^{-1}(\eta_{j+1})\rfloor+1,k\right) P_k x\right\|_n^\eta \\
	& \leq K N \left(\frac{\eta_{n+1}}{\eta_{j+2}}\right)^{-\nu}\left\|\Phi_A\left(\lfloor\widetilde\mu^{-1}(\eta_{j+1})\rfloor+1,k\right) P_kx\right\|_{j+2}^\eta \\
	& \le K^2 N \left(\frac{\eta_{n+1}}{\eta_{j+2}}\right)^{-\nu} \left\|x\right\|_k = K^2 N \left(\frac{\mu_{m}}{\mu_{k}}\frac{\eta_{n+1}}{\mu_m}\frac{\mu_k}{\eta_j}\frac{\eta_j}{\eta_{j+2}}\right)^{-\nu} \left\|x\right\|_k\\
	& \le K^2 N \theta^{2\nu} \left(\frac{\mu_{m}}{\mu_{k}}\right)^{-\nu} \left\|x\right\|_k,
\end{split}
\]
where we also used that $\frac{\mu_k}{\eta_j}\ge 1$ (as $k\ge \widetilde{\mu}^{-1}(\eta_j)$), $\frac{\eta_{n+1}}{\mu_m}\ge 1$ (since $\widetilde{\mu}^{-1}(\eta_{n+1})\ge m$) and $\frac{\eta_j}{\eta_{j+2}}\ge \theta^{-2}$ (see~\eqref{eq:bound2}).

\underline{Case II} - There is $n \in \N_0$ such that $$\lfloor\widetilde\mu^{-1}(\eta_{n+1})\rfloor +1> m \ge k
\ge \lfloor\widetilde\mu^{-1}(\eta_n)\rfloor+1.$$
By~\eqref{eq:bound-rates-2} 
and~\eqref{eq:ordinary-dich1-iii->i}, we get 
\[
\begin{split}
	\left\|\Phi_A(m, k) P_k x\right\|_m	\le K\|x\|_k 
	&=K\left (\frac{\mu_m}{\mu_k}\right )^{-\nu} \left (\frac{\mu_m}{\mu_k}\right )^{\nu} \|x\|_k \\
&\le K\left (\frac{\mu_m}{\mu_k}\right )^{-\nu} \left (\frac{\eta_{n+1}}{\eta_n}\right )^{\nu}\|x\|_k\\
& \le K \theta^{\nu} \left(\frac{\mu_{m}}{\mu_{k}}\right)^{-\nu} \|x\|_k.
\end{split}
\]
We conclude that~\eqref{eq:claim1} holds and the claim follows. 

Our next claim is that there exists $N_2\ge 1$ such that for $x \in X$ and $1\le m\le k$, we have
\begin{equation}\label{eq:claim2}
	\left\|\Phi_A(m, k) (\Id-P_k) x\right\|_m
	\leq N_2 \left(\frac{\mu_{k}}{\mu_{m}}\right)^{-\nu} \left\|x\right\|_k.
\end{equation}
We again split the proof of~\eqref{eq:claim2} into two cases.

\underline{Case I} - There are $n, j\in \N_0$ with $n<j$ such that 
\[
\lfloor \tilde \mu^{-1}(\eta_n)\rfloor +1\le m < \lfloor \tilde \mu^{-1}(\eta_{n+1})\rfloor +1 \le  \lfloor \tilde \mu^{-1}(\eta_{j})\rfloor +1\le k< \lfloor \tilde{\mu}^{-1}(\eta_{j+1})\rfloor +1.
\]
By~\eqref{eq:ordinary-dich2-iii->i} and~\eqref{eq:eta-dich2-iii->i} we have that 
\[
\begin{split}
& \quad \left\|\Phi_A(m, k) (\Id-P_k)x\right\|_m\\ 
	& =\left\|\Phi_A\left(m, \lfloor\widetilde\mu^{-1}(\eta_{n+1})\rfloor+1\right) (\Id- P_{\lfloor\widetilde\mu^{-1}(\eta_{n+1})\rfloor+1})\Phi_A\left(\lfloor\widetilde\mu^{-1}(\eta_{n+1})\rfloor +1, k\right) (\Id-P_k)x\right\|_m \\
 &\le K\left \|\Phi_A\left(\lfloor\widetilde\mu^{-1}(\eta_{n+1})\rfloor+1, \lfloor\widetilde\mu^{-1}(\eta_{j})\rfloor +1\right)
	\Phi_A\left(\lfloor\widetilde\mu^{-1}(\eta_{j})\rfloor+1,k\right) (\Id-P_k)x\right \|_{\lfloor \tilde {\mu}^{-1}(\eta_{n+1})\rfloor+1}\\
 &=K\left\|\Phi_{Q^{\mu,\eta}}\left(n+2,j+1\right) \Phi_A\left(\lfloor\widetilde\mu^{-1}(\eta_{j})\rfloor+1,k\right) (\Id-P_k)x\right\|_{n+2}^\eta\\
&\le KN \left (\frac{\eta_{j+1}}{\eta_{n+2}}\right )^{-\nu}\left \|\Phi_A\left(\lfloor\widetilde\mu^{-1}(\eta_{j})\rfloor+1,k\right) (\Id-P_k)x\right\|_{j+1}^\eta \\
&\le K^2 N \left (\frac{\eta_{j+1}}{\eta_{n+2}}\right )^{-\nu}\|x\|_k=K^2 N\left (\frac{\mu_k}{\mu_m}\right)^{-\nu}\left (\frac{\mu_m}{\mu_k} \frac{\eta_{j+1}}{\eta_{n+2}}\right)^{-\nu}\|x\|_k \\
&\le K^2N \theta^{2\nu}\left (\frac{\mu_k}{\mu_m}\right)^{-\nu}\|x\|_k,
\end{split}
\]
since  (see~\eqref{eq:bound2})
\[
\frac{\mu_m}{\mu_k} \frac{\eta_{j+1}}{\eta_{n+2}}=\frac{\mu_m}{\eta_{n+2}} \ \frac{\eta_{j+1}}{\mu_k} \ge \frac{\eta_n}{\eta_{n+2}} \ge \theta^{-2}.
\]

\underline{Case II} - There is $n\in \N_0$ such that
\[
\lfloor \tilde{\mu}^{-1}(\eta_n)\rfloor +1 \le m\le k<\lfloor \tilde{\mu}^{-1}(\eta_{n+1})\rfloor+1.
\]
By~\eqref{eq:bound-rates-2}
and~\eqref{eq:ordinary-dich2-iii->i} we have that 
\[
\begin{split}
\left\|\Phi_A(m, k) (\Id-P_k)x\right\|_m \le K\|x\|_k &=K\left (\frac{\mu_k}{\mu_m}\right )^{-\nu} \left (\frac{\mu_k}{\mu_m}\right )^{\nu}\|x\|_k\\
&\le K\left (\frac{\mu_k}{\mu_m}\right )^{-\nu} \left (\frac{\eta_{n+1}}{\eta_n}\right )^{\nu}\|x\|_k  \\
&\le K\theta^{\nu} \left (\frac{\mu_k}{\mu_m}\right )^{-\nu}\|x\|_k.
\end{split}
\]
We conclude that~\eqref{eq:claim2} holds and the claim follows.

We conclude that system~\eqref{eq:Lin} admits a $\mu$-dichotomy with respect to the sequence of norms $\left\{\|\cdot\|_k\right\}_{k \in \mathbb{N}}$ relative to the projections $\{P_k\}_{k \in \N}$. The proof of the theorem is completed.
\end{proof}

\begin{remark}
In the particular case where $\mu$ and $\eta$ are of the form $\mu_n=n+1$ and $\eta_n=h^n$, where $h\in \N\setminus \{1\}$, the version of Theorem~\ref{teo:main} was established in~\cite[Theorem 3.1]{DSS1}. Observe that in this case the $\eta$-dichotomy with respect to a sequence of norms is just an exponential dichotomy with respect to a sequence of norms.
\end{remark}

The following example shows that we cannot eliminate the assumption of an ordinary dichotomy in conditions~\eqref{teo:main-2} and~\eqref{teo:main-3} in the statement of Theorem~\ref{teo:main}.
\begin{example}
Take $X=\mathbb R$ and $\|\cdot \|_n=|\cdot |$ for $n\in \N$. Moreover, let $\mu_n=1+n$ and $\eta_n=2^n$ for $n\in \N_0$. Finally, we set
\[
A_n=\begin{cases}
n, & n=2^k-1 \ \text{for some $k\in \N\setminus \{1\}$;}\\
0, & \text{otherwise.}
\end{cases}
\]
Observe that 
\[
Q_n^{\mu, \eta}=\Phi_A(2^{n}, 2^{n-1})=0 \quad \text{for $n\in \N$,}
\]
yielding that~\eqref{eq:Lin-Q} admits an exponential dichotomy with respect to the sequence of norms $\left\{\|\cdot\|_k^{\eta}\right\}_{k \in \mathbb{N}}$. On the other hand, since $\sup_{n\in \N}A_n=+\infty$, it is easy to show that~\eqref{eq:Lin} does not admit $\mu$-dichotomy with respect to the sequence of norms $\left \{\| \cdot \|_k \right\}_{k\in \N}$. 
\end{example}

The following result shows that we can formulate additional conditions under which we can eliminate the assumption of an ordinary dichotomy in conditions~\eqref{teo:main-2} and~\eqref{teo:main-3} in the statement of Theorem~\ref{teo:main}.
\begin{theorem}\label{newthm}
Let $\mu$ and $\eta$ be growth rates with the property that there exists $\theta \ge 1$ such that~\eqref{eq:bound-rates-1-a} and~\eqref{eq:bound-rates-2} hold. Moreover, suppose that  $A_k$ is  an invertible operator for each $k\in \N$, and that there exists $K, a>0$ such that 
\begin{equation}\label{boundedgrowth1}
    \| \Phi_A(m, k)x\|_m \le K\left (\frac{\mu_m}{\mu_k}\right )^a \|x\|_k \quad \text{for $m\ge k$ and $x\in X$}
\end{equation}
and
\begin{equation}\label{boundedgrowth2}
    \| \Phi_A(m, k)x\|_m \le K\left (\frac{\mu_k}{\mu_m}\right )^a \|x\|_k \quad \text{for $m\le k$ and $x\in X$.}
\end{equation}
Then, there exists $K'>0$ such that 
\begin{equation}\label{boundedgrowth3}
  \|  \Phi_{Q^{\mu,\eta}}(m, k)x\|_m^\eta \le K' \left (\frac{\eta_m}{\eta_k}\right )^{a}\|x\|_k^\eta \quad \text{for $m\ge k$ and $x\in X$}
\end{equation}
and 
\begin{equation}\label{boundedgrowth4}
  \|  \Phi_{Q^{\mu,\eta}}(m, k)x\|_m^\eta \le K' \left (\frac{\eta_k}{\eta_m}\right )^{a}\|x\|_k^\eta \quad \text{for $m\le k$ and $x\in X$.}
\end{equation}
Moreover, the following properties are equivalent:
\begin{enumerate}
\item[(i)]  system~\eqref{eq:Lin} admits a $\mu$-dichotomy with respect to the sequence of norms $\left\{\|\cdot\|_k\right\}_{k \in \mathbb{N}}$;
\item[(ii)] system~\eqref{eq:Lin-Q} admits a $\eta$-dichotomy with respect to the norms $\left\{\|\cdot\|_k^{\eta}\right\}_{k \in \mathbb{N}}$.
\end{enumerate}
\end{theorem}
\begin{proof}
Take $m\ge k$ and $x\in X$. By~\eqref{eq:bound-rates-2}, \eqref{eq:bound2} and~\eqref{boundedgrowth1} we have that 
\[
\begin{split}
\|\Phi_{Q^{\mu,\eta}}(m, k)x\|_m^\eta 
&=\left \|\Phi_A \left (\lfloor\widetilde\mu^{-1}(\eta_{m-1})\rfloor +1,\,\lfloor\widetilde\mu^{-1}(\eta_{k-1})\rfloor +1\right ) x \right \|_{\lfloor \tilde{\mu}^{-1}(\eta_{m-1})\rfloor+1}\\
&\le K\left (\frac{\mu_{\lfloor\widetilde\mu^{-1}(\eta_{m-1})\rfloor+1}}{\mu_{\lfloor\widetilde\mu^{-1}(\eta_{k-1})\rfloor+1 }}\right )^a \|x\|_k^\eta \\
&\le K\left (\frac{\widetilde{\mu}(\widetilde\mu^{-1}(\eta_{m-1})+1)}{\eta_{k-1}}\right )^a  \|x\|_k^\eta\\
&\le K\theta^{2a} \left (\frac{\eta_{m-1}}{\eta_{k-1}}\right )^a \|x\|_k^\eta \\
&\le K\theta^{3a}\left (\frac{\eta_{m}}{\eta_{k}}\right )^a \|x\|_k^\eta.
\end{split}
\]
We conclude that~\eqref{boundedgrowth3} holds (with $K'=K\theta^{3a}$). Similarly, one can establish~\eqref{boundedgrowth4}.

We now establish the equivalence between $(i)$ and $(ii)$. It follows from Theorem~\ref{teo:main} that it is sufficient to prove that $(ii)$ implies $(i)$. Assume that system~\eqref{eq:Lin-Q} admits a $\eta$-dichotomy with respect to the norms $\left\{\|\cdot\|_k^{\eta}\right\}_{k \in \mathbb{N}}$ and projections $\widetilde{P}_k$, $k\in \N$. Thus, there exist $L, \nu>0$ such that 
\begin{equation}\label{munu1}
    \|\Phi_{Q^{\mu, \eta}}(m, k)\tilde P_kx\|_m^\eta\le L\left(\frac{\eta_m}{\eta_k}\right )^{-\nu}\|x\|_k^\eta \quad \text{for $m\ge k$ and $x\in X$,}
\end{equation}
and 
\begin{equation}\label{munu2}
\|\Phi_{Q^{\mu, \eta}}(m, k)(\Id-\tilde P_k)x\|_m^\eta\le L\left(\frac{\eta_k}{\eta_m}\right )^{-\nu}\|x\|_k^\eta \quad \text{for $m\le k$ and $x\in X$.}
\end{equation}
For $k\in \N$ we define
\[
P_k:=\Phi_A(k, 1)\tilde P_1\Phi_A(1, k).
\]
Then, 
\[
\begin{split}
P_{k+1}A_k =\Phi_A(k+1, 1)\tilde P_1\Phi_A(1, k+1)A_k
&=\Phi_A(k+1, 1)\tilde P_1\Phi_A(1, k) \\
&=A_k\Phi_A(k, 1)\tilde P_1\Phi_A(1, k) \\
&=A_kP_k,
\end{split}
\]
for each $k\in \N$ yielding (od\ref{def:od1}). The invertibility of the operators $A_k$ immediately yields (od\ref{def:od3}). Clearly,
\[
P_{\lfloor \widetilde{\mu}^{-1}(\eta_{k-1})\rfloor+1}=\tilde P_k, \quad k\in \N.
\]
We claim that there exists $D>0$ such that 
\begin{equation}\label{Pk}
\|P_kx\|_k \le D\|x\|_k, \quad \text{for $k\in \N$ and $x\in X$.}
\end{equation}
To this end, we fix $k\in \N$ and choose $n\in \N_0$ such that 
\[
\lfloor \widetilde{\mu}^{-1}(\eta_{n+1})\rfloor+1 >k\ge \lfloor \widetilde{\mu}^{-1}(\eta_{n})\rfloor+1.
\]
Then,  using~\eqref{eq:bound-rates-1-a}, \eqref{eq:bound-rates-2},
\eqref{boundedgrowth1}, \eqref{boundedgrowth2} and~\eqref{munu1} (applied for $m=k=n+2$) we have that 
\[
\begin{split}
\|P_kx\|_k &= \left \|\Phi_A(k, \lfloor \widetilde{\mu}^{-1}(\eta_{n+1})\rfloor+1 )\widetilde{P}_{n+2}\Phi_A (\lfloor \widetilde{\mu}^{-1}(\eta_{n+1})\rfloor+1, k)x\right \|_k \\
&\le K \left (\frac{\mu_{\lfloor \widetilde{\mu}^{-1}(\eta_{n+1})\rfloor+1}}{\mu_k}\right )^a \left \|\widetilde{P}_{n+2}\Phi_A (\lfloor \widetilde{\mu}^{-1}(\eta_{n+1})\rfloor+1, k)x\right \|_{n+2}^\eta \\
&\le K\theta^a \left (\frac{\eta_{n+1}}{\mu_k}\right )^a\left \|\widetilde{P}_{n+2}\Phi_A (\lfloor \widetilde{\mu}^{-1}(\eta_{n+1})\rfloor+1, k)x\right \|_{n+2}^\eta \\
&\le K \theta^a \left(\frac{\eta_{n+1}}{\eta_n}\right )^a  \left \|\widetilde{P}_{n+2}\Phi_A (\lfloor \widetilde{\mu}^{-1}(\eta_{n+1})\rfloor+1, k)x\right \|_{n+2}^\eta \\
&\le K\theta^{2a}\left \|\widetilde{P}_{n+2}\Phi_A (\lfloor \widetilde{\mu}^{-1}(\eta_{n+1})\rfloor+1, k)x\right \|_{n+2}^\eta \\
&\le KL\theta^{2a}\left \|\Phi_A (\lfloor \widetilde{\mu}^{-1}(\eta_{n+1})\rfloor+1, k)x\right \|_{n+2}^\eta\\
&\le K^2L\theta^{2a} \left (\frac{\mu_{\lfloor \widetilde{\mu}^{-1}(\eta_{n+1})\rfloor+1}}{\mu_k}\right )^a\|x\|_k\\
&\le K^2L\theta^{3a} \left (
\frac{\mu_{\lfloor \widetilde{\mu}^{-1}(\eta_{n+1})\rfloor}}{\eta_n}\right )^a\|x\|_k\\
&\le K^2L\theta^{4a}\|x\|_k,
\end{split}
\]
yielding~\eqref{Pk}.

We now claim that there exists $L'>0$ such that
\begin{equation}\label{947}
\left \|\Phi_A(m, k)P_kx \right \|_m \le L'\left(\frac{\mu_m}{\mu_k}\right )^{-\nu}\|x\|_k \quad \text{for $m\ge k$ and $x\in X$.}
\end{equation}

We distinguish two cases.

\underline{Case I} - There are $n,j \in \N_0$, $j<n$ such that 
$$\lfloor\widetilde\mu^{-1}(\eta_{n+1})\rfloor +1> m \ge \lfloor\widetilde\mu^{-1}(\eta_n)\rfloor +1 \ge \lfloor\widetilde\mu^{-1}(\eta_{j+1})\rfloor +1>  k
\ge \lfloor\widetilde\mu^{-1}(\eta_j)\rfloor +1.$$
By~\eqref{eq:bound-rates-2}, \eqref{eq:bound2}, \eqref{boundedgrowth1} and~\eqref{munu1} we have that  
\[
\begin{split}
	& \quad \left\|\Phi_A(m, k) P_k x\right\|_m\\ 
	& =\left\|
	\Phi_A\left(m, \lfloor\widetilde\mu^{-1}(\eta_n)\rfloor +1\right) \Phi_A\left(\lfloor\widetilde\mu^{-1}(\eta_n)\rfloor +1, \lfloor\widetilde\mu^{-1}(\eta_{j+1})\rfloor +1\right)
	\Phi_A\left(\lfloor\widetilde\mu^{-1}(\eta_{j+1})\rfloor+1,k\right) P_k x\right\|_m \\
	& \le K \left (\frac{\mu_m}{\mu_{\lfloor \widetilde{\mu}^{-1}(\eta_n)\rfloor+1}}\right)^a
 \left\| \Phi_A\left(\lfloor\widetilde\mu^{-1}(\eta_n)\rfloor+1, \lfloor\widetilde\mu^{-1}(\eta_{j+1})\rfloor +1\right)
	\Phi_A\left(\lfloor\widetilde\mu^{-1}(\eta_{j+1})\rfloor +1,k\right) P_k x\right\|_{\lfloor\widetilde\mu^{-1}(\eta_n)\rfloor +1}\\
 & \le K \left (\frac{\eta_{n+1}}{\eta_n}\right)^a \left\|\Phi_{Q^{\mu,\eta}}\left(n+1,j+2\right)\tilde P_{j+2} \Phi_A\left(\lfloor\widetilde\mu^{-1}(\eta_{j+1})\rfloor+1,k\right)  x\right\|_{n+1}^\eta
 \\
	 &\leq K\theta^{a}\left\|\Phi_{Q^{\mu,\eta}}\left(n+1,j+2\right)\tilde P_{j+2} \Phi_A\left(\lfloor\widetilde\mu^{-1}(\eta_{j+1})\rfloor+1,k\right)  x\right\|_{n+1}^\eta\\
	& \leq K L\theta^{a}  \left(\frac{\eta_{n+1}}{\eta_{j+2}}\right)^{-\nu}\left\|\Phi_A\left(\lfloor\widetilde\mu^{-1}(\eta_{j+1})\rfloor+1,k\right) x\right\|_{j+2}^\eta \\
 & \leq K^2 L \theta^{a} \left(\frac{\eta_{n+1}}{\eta_{j+2}}\right)^{-\nu} \left(\frac{\mu_{\lfloor \widetilde{\mu}^{-1}(\eta_{j+1})\rfloor+1}}{\mu_k}\right )^a\|x\|_k\\
 & \leq K^2 L \theta^{a} \left(\frac{\eta_{n+1}}{\eta_{j+2}}\right)^{-\nu} \left(\frac{\mu_{\lfloor \widetilde{\mu}^{-1}(\eta_{j+1})\rfloor+1}}{\eta_j}\right )^a\|x\|_k\\
 &\leq K^2 L\theta^{4a} \left(\frac{\eta_{n+1}}{\eta_{j+2}}\right)^{-\nu}\|x\|_k\\
 &=K^2L\theta^{4a}\left (  \frac{\mu_m}{\mu_k} \frac{\eta_{n+1}}{\mu_m} \frac{\mu_k}{\eta_j}\frac{\eta_j}{\eta_{j+2}}\right)^{-\nu}\|x\|_k                            \\
 &\le K^2L\theta^{4a+2\nu} \left (\frac{\mu_m}{\mu_k}\right )^{-\nu}\|x\|_k,
\end{split}
\]
where in the last step we used that $\frac{\eta_{n+1}}{\mu_m}\ge 1$, $\frac{\mu_k}{\eta_j}\ge 1$ and $\frac{\eta_j}{\eta_{j+2}}\ge \theta^{-2}$.

\underline{Case II} - There is $n\in \N_0$ such that 
\[
\lfloor\widetilde\mu^{-1}(\eta_{n+1})\rfloor+1 > m\ge k \ge \lfloor\widetilde\mu^{-1}(\eta_n)\rfloor+1.
\]
It follows from~\eqref{eq:bound-rates-2}, \eqref{boundedgrowth1}
and~\eqref{Pk} that 
\[
\begin{split}
\left\|\Phi_A(m, k) P_k x\right\|_m \le K\left (\frac{\mu_m}{\mu_k}\right )^a\|P_k x\|_k &\le KD\left (\frac{\mu_m}{\mu_k}\right )^a\|x\|_k \\
&\le KD\left (\frac{\eta_{n+1}}{\eta_n} \right )^a \|x\|_k \\
&\le KD\theta^{a}\|x\|_k \\
&= KD\theta^{a}\left (\frac{\mu_m}{\mu_k}\right )^{-\nu}\left (\frac{\mu_m}{\mu_k} \right )^\nu\|x\|_k\\
&\le KD\theta^{a+\nu}\left (\frac{\mu_m}{\mu_k}\right )^{-\nu}\|x\|_k.
\end{split}
\]
We conclude that~\eqref{947} holds. Similarly, one can show that there exists $L''>0$ such that 
\[
\left \|\Phi_A(m, k)(\Id-P_k)x \right \|_m \le L'' \left (\frac{\mu_k}{\mu_m}\right )^{-\nu}\|x\|_k, \quad \text{for $m\le k$ and $x\in X$.}
\]
Consequently, the system~\eqref{eq:Lin} admits a $\mu$-dichotomy with respect to the sequence of norms $\left\{\|\cdot\|_k\right\}_{k \in \mathbb{N}}$ and the proof of the theorem is completed.
\end{proof}

\section{$\mu$-dichotomy}

We start with the following definition.
\begin{definition}
Let $\mu=(\mu_n)_{n\in \N_0}$ be a discrete growth rate.
We say that~\eqref{eq:Lin} (or the sequence $\{A_n\}_{n\in \N}$)  admits:
\begin{enumerate}
\item ordinary dichotomy if it admits ordinary dichotomy with respect to the sequence of norms $\left \{ \| \cdot \|_k \right \}_{k\in \N}$ where $\|\cdot \|_k=\|\cdot \|$ for each $k\in \N$;
\item $\mu$-dichotomy if it admits $\mu$-dichotomy with respect to the sequence of norms $\left \{ \| \cdot \|_k \right \}_{k\in \N}$ where $\|\cdot \|_k=\|\cdot \|$ for each $k\in \N$.
\end{enumerate}
\end{definition}

Letting $\|\cdot\|_k=\|\cdot\|$ in Theorem~\ref{teo:main}, for all $k \in \N$, we obtain the following immediate corollary.

\begin{corollary}\label{corollary:main2}
	The following assertions are equivalent:
	\begin{enumerate}[(i)] 	
			\item\label{cor:main2-1} $\mu$ is a growth rate satisfying 
			\begin{equation}\label{eq:bound-rates-1a}
				\frac{\mu_{n+1}}{\mu_n}\le \theta, \quad \text{for all} \quad n \in \N_0,
			\end{equation} 
			for some $\theta\ge 1$, and system~\eqref{eq:Lin} admits a $\mu$-dichotomy with respect to a sequence of projections $\{P_k\}_{k \in \N}$;
		\item\label{cor:main2-2} system~\eqref{eq:Lin} admits an ordinary dichotomy relative to a sequence of projections $\{P_k\}_{k \in \mathbb{N}}$ and, for any growth rate $\eta$ satisfying 
		\begin{equation}\label{eq:bound-rates-2a}
			\frac{\eta_{n+1}}{\eta_n} \le \theta, \quad \text{for all} \quad n \in \N_0,
		\end{equation}
		for some $\theta\ge 1$, system~\eqref{eq:Lin-Q} admits a $\eta$-dichotomy relative to the sequence of projections $\left\{P^\eta_k\right\}_{k \in \mathbb{N}}$, where $P^\eta_k=P_{\lfloor \widetilde{\mu}^{-1}(\eta_{k-1})\rfloor+1}$, for all $k \in \mathbb{N}$;
		\item\label{cor:main2-3} system~\eqref{eq:Lin} admits an ordinary dichotomy relative to a sequence of projections $\{P_k\}_{k \in \mathbb{N}}$ and there is a growth rate $\eta$ satisfying~\eqref{eq:bound-rates-2a} such that system~\eqref{eq:Lin-Q} admits a $\eta$-dichotomy relative to the sequence of projections $\left\{P^\eta_k\right\}_{k \in \mathbb{N}}$, where $P^\eta_k=P_{\lfloor \widetilde{\mu}^{-1}(\eta_{k-1})\rfloor+1}$, for all $k \in \mathbb{N}$.
	\end{enumerate} 
\end{corollary}
The following result follows readily from Theorem~\ref{newthm}. 
\begin{corollary}\label{cor-442}
Let $\mu$ and $\eta$ be growth rates with the property that there exists $\theta \ge 1$ such that~\eqref{eq:bound-rates-1-a} and~\eqref{eq:bound-rates-2} hold. Moreover, suppose that  $A_k$ is  an invertible operator for each $k\in \N$, and that there exists $K, a>0$ such that 
\[
    \| \Phi_A(m, k)\| \le K\left (\frac{\mu_m}{\mu_k}\right )^a \quad \text{for $m\ge k$}
\]
and
\[
    \| \Phi_A(m, k)\| \le K\left (\frac{\mu_k}{\mu_m}\right )^a \quad \text{for $m\le k$.}
\]
Then, there exists $K'>0$ such that 
\[
  \|  \Phi_{Q^{\mu,\eta}}(m, k)x\| \le K' \left (\frac{\eta_m}{\eta_k}\right )^{a} \quad \text{for $m\ge k$}
\]
and 
\[
  \|  \Phi_{Q^{\mu,\eta}}(m, k)\| \le K' \left (\frac{\eta_k}{\eta_m}\right )^{a} \quad \text{for $m\le k$.}
\]
Moreover, the following properties are equivalent:
\begin{enumerate}
\item[(i)]  system~\eqref{eq:Lin} admits a $\mu$-dichotomy;
\item[(ii)] system~\eqref{eq:Lin-Q} admits a $\eta$-dichotomy.
\end{enumerate}
\end{corollary}

\section{Dichotomy spectrum}
In this section we derive some consequences of Theorem~\ref{newthm} for the study of the Sacker-Sell spectrum and its generalizations. 

Throughout this section, we consider a growth rate $\mu$ with the property that~\eqref{eq:bound-rates-1-a} holds with some $\theta \ge 1$.  Moreover, we fix a sequence $\mathbb A=\{A_n\}_{n\in \N}\subset \mathcal B(X)$ of invertible operators. Finally, let $\eta_n=e^n$, $n\in \N_0$.

By $\Sigma_{\mu D, \mathbb A}$ we denote the set of all $\lambda \in \R$ with the property that the system
\begin{equation}\label{lin1}
    x_{n+1}=\left (\frac{\mu_{n+1}}{\mu_n}\right )^{-\lambda }A_n x_n  =: \tilde{A}_nx_n \quad n\in \N
\end{equation}
does not admit a $\mu$-dichotomy. 

By $\Sigma_{ED, \mathbb Q}$ we denote the set of all $\lambda \in \R$ with the property that the system
\begin{equation}\label{lin2}
    y_{n+1}=e^{-\lambda} Q_n^{\mu, \eta} y_n \quad n\in \N
\end{equation}
does not admit an exponential dichotomy.

The following is the main result of this section.
\begin{theorem}\label{dicspec}
Suppose that there exist $K, a>0$ such that
\begin{equation}\label{bg-1}
\|\Phi_A(m, k)\| \le K\left (\frac{\mu_m}{\mu_k}\right )^a \quad m\ge k,
\end{equation}
and 
\begin{equation}\label{bg-2}
\|\Phi_A(m, k)\| \le K\left (\frac{\mu_k}{\mu_m}\right )^a \quad m\le k.
\end{equation}
 Then, 
\[
\Sigma_{\mu D, \mathbb A}=\Sigma_{ED, \mathbb Q}.
\]
     
\end{theorem}

\begin{proof}
Take $\lambda \notin \Sigma_{\mu D, \mathbb A}$. Then \eqref{lin1} admits a $\mu$-dichotomy.  By Corollary~\ref{cor-442} we conclude that system \[ y_{n+1}=\left (\frac{\mu_{\lfloor \widetilde{\mu}^{-1}(e^n)\rfloor+1}}{\mu_{\lfloor \widetilde{\mu}^{-1}(e^{n-1})\rfloor+1}}\right )^{-\lambda}Q_n^{\mu, \eta} y_n =: \tilde Q_n^{\mu, \eta} y_n \quad n\in \N \] admits an exponential dichotomy. Observe that  
$$\tilde{Q}_n^{\mu,\eta}=\left (\frac{\mu_{\lfloor \widetilde{\mu}^{-1}(e^n)\rfloor+1}}{\mu_{\lfloor \widetilde{\mu}^{-1}(e^{n-1})\rfloor+1}}\right )^{-\lambda}\Phi_{A}(\lfloor \widetilde{\mu}^{-1}(e^n)\rfloor+1,\lfloor \widetilde{\mu}^{-1}(e^{n-1})\rfloor+1), \quad n\in \N.$$ Hence, there exists a sequence of projections $P_k$, $k\in \N$ on $X$ satisfying 
\[
P_{k+1}Q_k^{\mu, \eta}=Q_k^{\mu, \eta} P_k \quad k\in \N,
\]
and constants $K, \nu>0$ such that 
\begin{equation}\label{8:52}
\left (\frac{\mu_{\lfloor \widetilde{\mu}^{-1}(e^{m-1})\rfloor+1}}{\mu_{\lfloor \widetilde{\mu}^{-1}(e^{k-1})\rfloor+1}}\right )^{-\lambda}\left \| \Phi_{Q^{\mu, \eta}}(m, k)P_k \right \| \le Ke^{-\nu (m-k)} \quad m\ge  k\ge 1,
\end{equation}
and 
\begin{equation}\label{8:53}
\left (\frac{\mu_{\lfloor \widetilde{\mu}^{-1}(e^{m-1})\rfloor+1}}{\mu_{\lfloor \widetilde{\mu}^{-1}(e^{k-1})\rfloor+1}}\right )^{-\lambda}\left \| \Phi_{Q^{\mu, \eta}}(m, k)(\Id-P_k) \right \| \le Ke^{-\nu (k-m)} \quad 1\le m\le k.
\end{equation}
Next, note that 
\[
\left (\frac{\mu_{\lfloor \widetilde{\mu}^{-1}(e^{m-1})\rfloor+1}}{\mu_{\lfloor \widetilde{\mu}^{-1}(e^{k-1})\rfloor+1}}\right )^{-\lambda}
=e^{-\lambda (m-k)} 
\left (\frac{\widetilde{\mu} (\lfloor \widetilde{\mu}^{-1}(e^{m-1})\rfloor+1)}{\widetilde{\mu}(\widetilde{\mu}^{-1}(e^{m-1}))}\right )^{-\lambda} \left (\frac{\widetilde{\mu}(\widetilde{\mu}^{-1}(e^{k-1}))}{\widetilde{\mu} (\lfloor \widetilde{\mu}^{-1}(e^{k-1})\rfloor+1)}\right )^{-\lambda},
\]
for $m, k\in \N$.
Since 
\[
1\le \frac{\widetilde{\mu} (\lfloor \widetilde{\mu}^{-1}(e^{m-1})\rfloor+1)}{\widetilde{\mu}(\widetilde{\mu}^{-1}(e^{m-1}))}\le \theta^2
\]
and
\[
\theta^{-2} \le \frac{\widetilde{\mu}(\widetilde{\mu}^{-1}(e^{k-1}))}{\widetilde{\mu} (\lfloor \widetilde{\mu}^{-1}(e^{k-1})\rfloor+1)}\le 1,
\]
\eqref{8:52} and~\eqref{8:53} imply that there exists $C>0$ independent on $m$ and $k$ such that 
\[
e^{-\lambda (m-k)}\left \| \Phi_{Q^{\mu, \eta}}(m, k)P_k \right \| \le Ce^{-\nu (m-k)} \quad m\ge  k\ge 1,
\]
and
\[
e^{-\lambda (m-k)}\left \| \Phi_{Q^{\mu, \eta}}(m, k)(\Id-P_k) \right \| \le Ce^{-\nu (k-m)} \quad 1\le m\le k.
\]
We conclude that~\eqref{lin2} admits an exponential dichotomy. Hence, $\lambda \notin \Sigma_{ED, \mathbb Q}$ and, consequently, $\Sigma_{ED, \mathbb Q} \subset \Sigma_{\mu D, \mathbb A}$. Similarly, one can establish the reverse inclusion yielding the desired conclusion.
\end{proof}

\begin{remark}
Theorem~\ref{dicspec} enables us to deduce the structure of the spectrum $\Sigma_{\mu D, \mathbb A}$ by using known results about the Sacker-Sell spectrum. In fact, suppose that $X=\mathbb R^n$ and that there exist $K, a>0$ such that~\eqref{bg-1} and~\eqref{bg-2} hold. By Corollary~\ref{cor-442} we have
\[
\sup_{n\in \N}\|Q_n^{\mu, \eta}\|<+\infty \quad \text{and} \quad \sup_{n\in \N}\|(Q_n^{\mu, \eta})^{-1}\|<+\infty.
\]
By~\cite[Theorem 3.4]{AS} we see that $\Sigma_{ED, \mathbb Q}$ is a union of at most $n$ closed intervals in $\mathbb R$. Hence, $\Sigma_{\mu D, \mathbb A}$ has the same structure. 
\end{remark}

\begin{remark}
We note that the modifications of the Sacker-Sell spectrum associated with growth rates have been studied in~\cite{BDV, JG, Silva}. However, in these works, the spectrum is introduced with respect to the notion of a \emph{nonuniform} $\mu$-dichotomy. 
\end{remark}

\section{Nonautonomous topological linearization}

Throughout this section, we consider $\mu$ to be a growth rate that satisfies~\eqref{eq:bound-rates-1-a} with some $\theta \ge 1$. Besides~\eqref{eq:Lin} we will consider its nonlinear perturbations of the form
\[
x_{n+1}=A_n x_n+g_n(x_n) \quad n\in \N,
\]
where $g_n\colon X\to X$, $n\in \mathbb N$ is a sequence of maps.

\begin{theorem}\label{NTL}
Assume that the following conditions hold:
\begin{enumerate}
\item \eqref{eq:Lin} admits a $\mu$-dichotomy with respect to a sequence of norms $\left \{\|\cdot \|_m \right \}_{m\in \N}$, where each $\|\cdot \|_m$ is equivalent to $\| \cdot \|$. Moreover, suppose that the operators  $A_n$ are invertible;
\item there exist $K, a>0$ such that~\eqref{boundedgrowth1} and~\eqref{boundedgrowth2} hold;
\item there exists $M>0$ such that 
\begin{equation}\label{BND}
    \|g_n(x)\|_{n+1}\le M \frac{\mu_n'}{\mu_n} \quad \text{for $n\in \mathbb N$ and $x\in X$,}
\end{equation}
where $\mu_n'=\mu_{n+1}-\mu_n$;
\item  there exists $c>0$ such that 
\begin{equation}\label{LIP}
    \|g_n(x)-g_n(y)\|_{n+1}\le c\frac{\mu_n'}{\mu_n}\|x-y\|_n, \quad \text{for $n\in \N$ and $x,y\in X$.}
\end{equation}
\end{enumerate}
Then, provided that $c$ is sufficiently small, there exists a sequence $(h_n)_{n\in \N}$ of homeomorphisms $h_n\colon X \to X$ such that the following holds:
\begin{enumerate}
\item $h_{n+1}\circ (A_n+g_n)=A_n \circ h_n$ for $n\in \N$;
\item there exists $D>0$ such that for $x\in X$ and $n\in \N$,
\begin{equation}\label{Id+}
\|h_n(x)-x\|_n \le D \quad \text{and} \quad \|h_n^{-1}(x)-x\|_n \le D;
\end{equation}
\item there exist $\bar K, \varrho>0$ and $\delta>0$ such that 
\begin{equation}\label{Holdercont}
\|h_n(x)-h_n(y)\|_n \le \bar K\|x-y\|_n^\varrho \quad \text{and} \quad \|h_n^{-1}(x)-h_n^{-1}(y)\| \le \bar K\|x-y\|_n^{\varrho},
\end{equation}
for $n\in \N$ and $x, y\in X$ with $\|x-y\|_n\le \delta$.
\end{enumerate}
\end{theorem}

\begin{proof}
We first claim that provided $c$ is sufficiently small, $A_n+g_n$ is a homeomorphism for each $n\in \N$. To this end, let us fix an arbitrary $y\in X$. We claim that the map $x\mapsto A_n^{-1}y-A_n^{-1}g_n(x)$ is a contraction on $(X, \|\cdot \|_n)$ for each $n\in \N$. Indeed, observe that for $x_1, x_2\in X$ we have (using~\eqref{eq:bound-rates-1-a}, \eqref{boundedgrowth2} and~\eqref{LIP}) that
\[
\begin{split}
\|A_n^{-1}y-A_n^{-1}g_n(x_1)-(A_n^{-1}y-A_n^{-1}g_n(x_2))\|_n &=\|A_n^{-1}g_n(x_1)-A_n^{-1}g_n(x_2)\|_n \\
&\le K \left (\frac{\mu_{n+1}}{\mu_n}\right )^a \|g_n(x_1)-g_n(x_2)\|_{n+1} \\
&\le cK\left (\frac{\mu_{n+1}}{\mu_n}\right )^a \frac{\mu_{n}'}{\mu_n}\|x_1-x_2\|_n \\
&\le cK\left (\frac{\mu_{n+1}}{\mu_n}\right )^{a+1}\|x_1-x_2\|_n \\
&\le cK\theta^{a+1}\|x_1-x_2\|_n.
\end{split}
\]
Therefore, provided that $cK\theta^{a+1}<1$, we find that
$x\mapsto A_n^{-1}y-A_n^{-1}g_n(x)$ is a contraction, which produces the existence of a unique $x\in X$ satisfying $x=A_n^{-1}y-A_n^{-1}g_n(x)$, i.e. $A_n x+g_n(x)=y$.

For $m, n\in \mathbb N$ we set
\begin{equation}\label{gmn}
\mathcal G(m, n):=\begin{cases}
(A_{m-1}+g_{m-1})\circ \ldots \circ (A_n+g_n) & m>n; \\
\Id & m=n; \\
(A_m+g_m)^{-1}\circ \ldots \circ (A_{n-1}+g_{n-1})^{-1} & m<n.
\end{cases}
\end{equation}
We claim that there exist $\tilde K, \tilde a>0$ such that 
\begin{equation}\label{BND-1}
    \|\mathcal G(m, n)(x)-\mathcal G(m, n)(y)\|_m \le \tilde K\left (\frac{\mu_m}{\mu_n}\right )^{\tilde a}\|x-y\|_n, \quad \text{for $m\ge n$ and $x, y\in X$.}
\end{equation}
To  this end, we first note that 
\begin{equation}\label{gj}
\mathcal G(m, n)(x)=\Phi_A(m, n)x+\sum_{j=n}^{m-1} \Phi_A (m, j+1)g_j(\mathcal G(j, n)(x)),
\end{equation}
for $x\in X$ and $m\ge n \ge 1$.  Then, by~\eqref{boundedgrowth1} and~\eqref{BND} we have that 
\[
\begin{split}
&\|\mathcal G(m, n)(x)-\mathcal G(m, n)(y)\|_m\\
&\le K\left (\frac{\mu_m}{\mu_n}\right)^a\|x-y\|_n +Kc\sum_{j=n}^{m-1}\left (\frac{\mu_m}{\mu_{j+1}}\right )^a \frac{\mu_j'}{\mu_j}\|\mathcal G(j, n)(x)-\mathcal G(j, n)(y)\|_j \\
& \le K\left (\frac{\mu_m}{\mu_n}\right)^a \|x-y\|_n +  Kc\sum_{j=n}^{m-1}\left (\frac{\mu_m}{\mu_{j}}\right )^a \frac{\mu_j'}{\mu_j}\|\mathcal G(j, n)(x)-\mathcal G(j, n)(y)\|_j.
\end{split}
\]
yielding that 
\[
\begin{split}
&\left(\frac{\mu_n}{\mu_m}\right )^a\|\mathcal G(m, n)(x)-\mathcal G(m, n)(y)\|_m \\
&\le K\|x-y\|_n
+Kc\sum_{j=n}^{m-1}\left (\frac{\mu_n}{\mu_j}\right )^a \frac{\mu_j'}{\mu_j}\|\mathcal G(j, n)(x)-\mathcal G(j, n)(y)\|_j.
\end{split}
\]
Consequently, the discrete Gronwall's lemma gives that
\[
\left(\frac{\mu_n}{\mu_m}\right )^a\|\mathcal G(m, n)(x)-\mathcal G(m, n)(y)\|_m\le K e^{Kc\sum_{j=n}^{m-1}\frac{\mu_j'}{\mu_j}}\|x-y\|_n.
\]
Since (see~\cite[Lemma 3.1]{DSV})  \begin{equation}\label{nm-1}\sum_{j=n}^{m-1}\frac{\mu_j'}{\mu_j}\le \theta \log \left (\frac{\mu_m}{\mu_n}\right ),\end{equation} we conclude that 
\[
\|\mathcal G(m, n)(x)-\mathcal G(m, n)(y)\|_m \le K \left (\frac{\mu_m}{\mu_n}\right )^{a+cK\theta}\|x-y\|_n.
\]
We see that~\eqref{BND-1} holds with $\tilde K=K$ and $\tilde a:=a+cK\theta$. 

Similarly, using that 
\begin{equation}\label{Gmn}
\mathcal G(m, n)(x)=\Phi_A(m, n)(x)-\sum_{j=m}^{n-1}\Phi_A(m, j+1)g_j(\mathcal G(j, n)(x)) 
\end{equation}
for $1\le m\le n$ and $x\in X$, one can show that 
\begin{equation}\label{BND-2}
    \|\mathcal G(m, n)(x)-\mathcal G(m, n)(y)\|_m \le \tilde K\left (\frac{\mu_n}{\mu_m}\right )^{\tilde a}\|x-y\|_n, \quad \text{for $m\le n$ and $x, y\in X$.}
\end{equation}
Set $\eta_n=e^n$, $n\in \mathbb N_0$. Moreover, for $n\in \N$ we define $f_n\colon X\to X$ by
\begin{equation}\label{fnx}
f_n(x):=\sum_{j=\lfloor \widetilde{\mu}^{-1}(\eta_{n-1})\rfloor+1}^{\lfloor \widetilde{\mu}^{-1}(\eta_{n})\rfloor}\Phi_A(\lfloor \widetilde{\mu}^{-1}(\eta_{n})\rfloor+1, j+1)g_j (\mathcal G(j, \lfloor \widetilde{\mu}^{-1}(\eta_{n-1})\rfloor+1)(x)), 
\end{equation}
for $x\in X$.
We claim that there exists $C_1>0$ such that 
\begin{equation}\label{BND-3}
    \|f_n(x)\|_{\lfloor \widetilde{\mu}^{-1}(\eta_{n})\rfloor+1}\le C_1, \quad  \text{for $n\in \mathbb N$ and  $x\in X$.}
\end{equation}
Indeed, by~\eqref{eq:bound2}, \eqref{boundedgrowth1}, \eqref{BND} and~\eqref{nm-1} we have 
\[
\begin{split}
 \|f_n(x)\|_{\lfloor \widetilde{\mu}^{-1}(\eta_{n})\rfloor+1}&\le KM\sum_{j=\lfloor \widetilde{\mu}^{-1}(\eta_{n-1})\rfloor+1}^{\lfloor \widetilde{\mu}^{-1}(\eta_{n})\rfloor}\left (\frac{\mu_{\lfloor \widetilde{\mu}^{-1}(\eta_{n})\rfloor+1}}{\mu_{j+1}}\right )^a \frac{\mu_j'}{\mu_j} \\
 &\le KM\left (\frac{\mu_{\lfloor \widetilde{\mu}^{-1}(\eta_{n})\rfloor+1}}{\mu_{\lfloor \widetilde{\mu}^{-1}(\eta_{n-1})\rfloor +2}}\right )^a \sum_{j=\lfloor \widetilde{\mu}^{-1}(\eta_{n-1})\rfloor+1}^{\lfloor \widetilde{\mu}^{-1}(\eta_{n})\rfloor}\frac{\mu_j'}{\mu_j}\\
 &\le KM\theta \left (\frac{\mu_{\lfloor \widetilde{\mu}^{-1}(\eta_{n})\rfloor+1}}{\mu_{\lfloor \widetilde{\mu}^{-1}(\eta_{n-1})\rfloor +2}}\right )^a \log  \left (\frac{\mu_{\lfloor \widetilde{\mu}^{-1}(\eta_{n})\rfloor+1}}{\mu_{\lfloor \widetilde{\mu}^{-1}(\eta_{n-1})\rfloor+1}}\right ) \\
 &\le  KM\theta \left (\frac{\mu_{\lfloor \widetilde{\mu}^{-1}(\eta_{n})\rfloor+1}}{\mu_{\lfloor \widetilde{\mu}^{-1}(\eta_{n-1})\rfloor +1}}\right )^a \log  \left (\frac{\mu_{\lfloor \widetilde{\mu}^{-1}(\eta_{n})\rfloor+1}}{\mu_{\lfloor \widetilde{\mu}^{-1}(\eta_{n-1})\rfloor+1}}\right )\\
 &\le KM\theta (3\theta^2)^a \log (3 \theta^2),
\end{split}
\]
as
\begin{equation}\label{5-01}
\frac{\mu_{\lfloor \widetilde{\mu}^{-1}(\eta_{n})\rfloor+1}}{\mu_{\lfloor \widetilde{\mu}^{-1}(\eta_{n-1})\rfloor+1}}\le \frac{\widetilde{\mu}(\widetilde{\mu}^{-1}(\eta_n)+1)}{\eta_{n-1}}=\frac{\widetilde{\mu}(\widetilde{\mu}^{-1}(\eta_n)+1)}{\eta_n} \ \frac{\eta_{n}}{\eta_{n-1}}\le 3\theta^2.
\end{equation}  
We conclude that~\eqref{BND-3} holds.

Furthermore, we claim that there exists $C_2>0$ such that 
\begin{equation}\label{LIP2}
    \|f_n(x)-f_n(y)\|_{\lfloor \widetilde{\mu}^{-1}(\eta_{n})\rfloor+1}\le cC_2\|x-y\|_{\lfloor \widetilde{\mu}^{-1}(\eta_{n-1})\rfloor+1} \quad \text{for $x, y\in X$ and $n\in \mathbb N$.}
\end{equation}
By~\eqref{boundedgrowth1}, \eqref{LIP}, \eqref{BND-1} and~\eqref{5-01}
we have by denoting $x_j:= \mathcal G(j, \lfloor \widetilde{\mu}^{-1}(\eta_{n-1})\rfloor+1)(x)$ and $y_j:=\mathcal G(j, \lfloor \widetilde{\mu}^{-1}(\eta_{n-1})\rfloor+1)(y)$ that 
\begin{equation}\label{fnxy}
\begin{split}
&\|f_n(x)-f_n(y)\|_{\lfloor \widetilde{\mu}^{-1} (\eta_{n})\rfloor+1} \\
&\le  K\sum_{j=\lfloor \widetilde{\mu}^{-1}(\eta_{n-1})\rfloor+1}^{\lfloor \widetilde{\mu}^{-1}(\eta_{n})\rfloor}\left (\frac{\mu_{\lfloor \widetilde{\mu}^{-1}(\eta_{n})\rfloor+1}}{\mu_{j+1}}\right )^a \|g_j (x_j)-g_j (y_j)\|_{j+1} \\
&\le Kc\left (\frac{\mu_{\lfloor \widetilde{\mu}^{-1}(\eta_{n})\rfloor+1}}{\mu_{\lfloor \widetilde{\mu}^{-1}(\eta_{n-1})\rfloor+1}}\right )^a\sum_{j=\lfloor \widetilde{\mu}^{-1}(\eta_{n-1})\rfloor+1}^{\lfloor \widetilde{\mu}^{-1}(\eta_{n})\rfloor}\frac{\mu_j'}{\mu_j}\|x_j-y_j\|_j \\
&\le K\tilde K c\left (\frac{\mu_{\lfloor \widetilde{\mu}^{-1}(\eta_{n})\rfloor+1}}{\mu_{\lfloor \widetilde{\mu}^{-1}(\eta_{n-1})\rfloor+1}}\right )^a \|x-y\|_{\lfloor \widetilde{\mu}^{-1}(\eta_{n-1})\rfloor+1}
\sum_{j=\lfloor \widetilde{\mu}^{-1}(\eta_{n-1})\rfloor+1}^{\lfloor \widetilde{\mu}^{-1}(\eta_{n})\rfloor}\frac{\mu_j'}{\mu_j}\left ( \frac{\mu_j}{\mu_{\lfloor \widetilde{\mu}^{-1}(\eta_{n-1})\rfloor+1}}\right )^{\tilde a}\\
&\le K\tilde K c \theta \|x-y\|_{\lfloor \widetilde{\mu}^{-1}(\eta_{n-1})\rfloor+1} \left (\frac{\mu_{\lfloor \widetilde{\mu}^{-1}(\eta_{n})\rfloor+1}}{\mu_{\lfloor \widetilde{\mu}^{-1}(\eta_{n-1})\rfloor+1}}\right )^{a+\tilde a }\log \left (\frac{\mu_{\lfloor \widetilde{\mu}^{-1}(\eta_{n})\rfloor+1}}{\mu_{\lfloor \widetilde{\mu}^{-1}(\eta_{n-1})\rfloor+1}}\right ) \\
&\le K\tilde K c \theta (3\theta^2)^{a+\tilde a} \log (3\theta^2)  \|x-y\|_{\lfloor \widetilde{\mu}^{-1}(\eta_{n-1})\rfloor+1},
\end{split}
\end{equation}
which yields~\eqref{LIP2}.

Let $Q_n^{\mu, \eta}$ be given by~\eqref{eq:Q-mu-eta}. By Theorem~\ref{newthm} we find that~\eqref{eq:Lin-Q} admits an exponential dichotomy with respect to the sequence of norms $\left \{\|\cdot \|_n^\eta \right \}_{n\in \N}$. Moreover, \eqref{boundedgrowth3} and~\eqref{boundedgrowth4} hold with some $K'>0$.  It follows from Corollary~\ref{BVcor} that there is a sequence of homeomorphisms $\psi_n\colon X\to X$, $n\in \mathbb N$ such that 
\begin{equation}\label{lin}
    \psi_{n+1}\circ \mathcal G(\lfloor \widetilde{\mu}^{-1}(\eta_{n})\rfloor+1, \lfloor \widetilde{\mu}^{-1}(\eta_{n-1})\rfloor+1)=\Phi_A (\lfloor \widetilde{\mu}^{-1}(\eta_{n})\rfloor+1, \lfloor \widetilde{\mu}^{-1}(\eta_{n-1})\rfloor+1)\circ \psi_n,
\end{equation}
for $\mathbb N$.  Moreover, there exist $D, \varrho>0$ such that 
\begin{equation}\label{holder1}
\|\psi_n(x)-\psi_n(y)\|_{\lfloor \widetilde{\mu}^{-1}(\eta_{n-1})\rfloor+1}\le D\|x-y\|_{\lfloor \widetilde{\mu}^{-1}(\eta_{n-1})\rfloor+1}^\varrho
\end{equation}
and
\begin{equation}\label{holder2}
\|\psi_n^{-1}(x)-\psi_n^{-1}(y)\|_{\lfloor \widetilde{\mu}^{-1}(\eta_{n-1})\rfloor+1}\le D\|x-y\|_{\lfloor \widetilde{\mu}^{-1}(\eta_{n-1})\rfloor+1}^\varrho,
\end{equation}
for $n\in \mathbb N$ and $x, y\in X$ such that $\|x-y\|_{\lfloor \widetilde{\mu}^{-1}(\eta_{n-1})\rfloor+1}\le 1$. Finally, 
\begin{equation}\label{6:10}
\|\psi_n(x)-x\|_{\lfloor \widetilde{\mu}^{-1}(\eta_{n-1})\rfloor+1}\le D \quad \text{and} \quad \|\psi_n^{-1}(x)-x\|_{\lfloor \widetilde{\mu}^{-1}(\eta_{n-1})\rfloor+1}\le D, 
\end{equation}
for $n\in \N$ and $x\in X$. Take $k\in \N$, choose   $n\in \N_0$ such that 
\begin{equation}\label{nn+1}\lfloor\widetilde\mu^{-1}(\eta_{n+1})\rfloor +1>  k
\ge \lfloor\widetilde\mu^{-1}(\eta_n)\rfloor+1,\end{equation}
and set
\begin{equation}\label{hk}
h_k:=\Phi_A (k, \lfloor\widetilde\mu^{-1}(\eta_{n})\rfloor+1)\circ \psi_{n+1}\circ \mathcal G(\lfloor\widetilde\mu^{-1}(\eta_{n})\rfloor+1, k).
\end{equation}
Clearly, $h_k$ is a homeomorphism for each $k\in \N$. We claim that 
\begin{equation}\label{linearization}
    h_{k+1}\circ (A_k+g_k)=A_k \circ h_k, \quad \text{for $k\in \N$.}
\end{equation}
To this end, let us fix $k\in \N$ and choose $n\in \N_0$ so that~\eqref{nn+1} holds.   If $\lfloor\widetilde\mu^{-1}(\eta_{n+1})\rfloor +1>k+1$, then 
\[
\begin{split}
&h_{k+1}\circ (A_k+g_k) \\
&=\Phi_A (k+1, \lfloor\widetilde\mu^{-1}(\eta_{n})\rfloor+1)\circ \psi_{n+1}\circ \mathcal G(\lfloor\widetilde\mu^{-1}(\eta_{n})\rfloor+1, k+1)\circ \mathcal G(k+1, k)\\
&=A_k \circ \Phi_A (k, \lfloor\widetilde\mu^{-1}(\eta_{n})\rfloor+1)\circ \psi_{n+1}\circ \mathcal G(\lfloor\widetilde\mu^{-1}(\eta_{n})\rfloor+1, k) \\
&=A_k \circ h_k,
\end{split}
\]
yielding~\eqref{linearization}.  On the other hand, if $\lfloor\widetilde\mu^{-1}(\eta_{n+1})\rfloor +1=k+1$, then using~\eqref{lin} we get that 
\[
\begin{split}
&h_{k+1}\circ (A_k+g_k) \\
&=\psi_{n+2} \circ \mathcal G(\lfloor\widetilde\mu^{-1}(\eta_{n+1})\rfloor +1, \lfloor\widetilde\mu^{-1}(\eta_{n+1})\rfloor) \\
&=\psi_{n+2}\circ \mathcal G(\lfloor\widetilde\mu^{-1}(\eta_{n+1})\rfloor +1, \lfloor\widetilde\mu^{-1}(\eta_{n})\rfloor +1)\circ \mathcal G(\lfloor\widetilde\mu^{-1}(\eta_{n})\rfloor +1, \lfloor\widetilde\mu^{-1}(\eta_{n+1})\rfloor)\\
&=\Phi_A (\lfloor\widetilde\mu^{-1}(\eta_{n+1})\rfloor +1, \lfloor\widetilde\mu^{-1}(\eta_{n})\rfloor +1)\circ \psi_{n+1}\circ \mathcal G(\lfloor\widetilde\mu^{-1}(\eta_{n})\rfloor +1, \lfloor\widetilde\mu^{-1}(\eta_{n+1})\rfloor)\\
&=A_k\circ \Phi_A (\lfloor\widetilde\mu^{-1}(\eta_{n+1})\rfloor , \lfloor\widetilde\mu^{-1}(\eta_{n})\rfloor +1)\circ \psi_{n+1}\circ \mathcal G(\lfloor\widetilde\mu^{-1}(\eta_{n})\rfloor +1, \lfloor\widetilde\mu^{-1}(\eta_{n+1})\rfloor)\\
&=A_k\circ h_k.
\end{split}
\]
Thus~\eqref{linearization} holds in this case as well.

Take $k\in \N$ and choose $n\in \N_0$ so that~\eqref{nn+1} holds.  By~\eqref{boundedgrowth1}, \eqref{BND-2} and~\eqref{holder1} we have (writing $x'=\mathcal G(\lfloor\widetilde\mu^{-1}(\eta_{n})\rfloor+1, k)(x)$ and $y'=\mathcal G(\lfloor\widetilde\mu^{-1}(\eta_{n})\rfloor+1, k)(y)$) that 
\[
\begin{split}
\|h_k(x)-h_k(y)\|_k 
&\le K \left (\frac{\mu_k}{\mu_{\lfloor\widetilde\mu^{-1}(\eta_{n})\rfloor+1}}\right )^a \left \| \psi_{n+1}( x')- \psi_{n+1}( y')\right \|_{\lfloor \widetilde{\mu}^{-1}(\eta_{n})\rfloor+1} \\
&\le KD \left (\frac{\mu_{\lfloor\widetilde\mu^{-1}(\eta_{n+1})\rfloor}}{\mu_{\lfloor\widetilde\mu^{-1}(\eta_{n})\rfloor+1}}\right )^a \|x'-y'\|_{\lfloor \widetilde{\mu}^{-1}(\eta_{n})\rfloor+1}^\varrho \\
&\le KD \tilde K^\varrho\left (\frac{\mu_{\lfloor\widetilde\mu^{-1}(\eta_{n+1})\rfloor}}{\mu_{\lfloor\widetilde\mu^{-1}(\eta_{n})\rfloor+1}}\right )^a \left (\frac{\mu_k}{\mu_{\lfloor \widetilde{\mu}^{-1}(\eta_{n})\rfloor+1}}\right )^{\varrho \tilde a}  \|x-y\|_k^\varrho \\
&\le KD \tilde K^\varrho\left (\frac{\mu_{\lfloor\widetilde\mu^{-1}(\eta_{n+1})\rfloor}}{\mu_{\lfloor\widetilde\mu^{-1}(\eta_{n})\rfloor+1}}\right )^a \left (\frac{\mu_{\lfloor\widetilde\mu^{-1}(\eta_{n+1})\rfloor}}{\mu_{\lfloor \widetilde{\mu}^{-1}(\eta_{n})\rfloor+1}}\right )^{\varrho \tilde a}  \|x-y\|_k^\varrho \\
&\le KD\tilde K^\varrho e^{a+\varrho \tilde a}\|x-y\|_k^\varrho,
\end{split}
\]
for $x, y\in X$ such that $\|x-y\|_k \le \frac{1}{\tilde Ke^{\tilde a}}$ which ensures that $\|x'-y'\|_{\lfloor \widetilde{\mu}^{-1}(\eta_{n})\rfloor+1}\le 1$. This yields the first estimate in~\eqref{Holdercont}. Similarly one can establish the second.

Next,  we again fix $k\in \N$ and choose $n\in \N_0$ satisfying~\eqref{nn+1}. Then, 
writing $\bar{\psi}_n:=\psi_n-\Id$ we have 
\[
\begin{split}
h_k-\Id &=\Phi_A (k, \lfloor\widetilde\mu^{-1}(\eta_{n})\rfloor+1)\circ \psi_{n+1}\circ \mathcal G(\lfloor\widetilde\mu^{-1}(\eta_{n})\rfloor+1, k)-\Id \\
&=\Phi_A (k, \lfloor\widetilde\mu^{-1}(\eta_{n})\rfloor+1)\circ \bar \psi_{n+1}\circ \mathcal G(\lfloor\widetilde\mu^{-1}(\eta_{n})\rfloor+1, k)-\Id \\
&\phantom{=}+\Phi_A (k, \lfloor\widetilde\mu^{-1}(\eta_{n})\rfloor+1)\circ  \mathcal G(\lfloor\widetilde\mu^{-1}(\eta_{n})\rfloor+1, k) \\
&=\Phi_A (k, \lfloor\widetilde\mu^{-1}(\eta_{n})\rfloor+1)\circ \bar \psi_{n+1}\circ \mathcal G(\lfloor\widetilde\mu^{-1}(\eta_{n})\rfloor+1, k)\\
&\phantom{=}-\sum_{j=\lfloor\widetilde\mu^{-1}(\eta_{n})\rfloor+1}^{k-1}\Phi_A(k, j+1)g_j(\mathcal G(j, k)(\cdot)),
\end{split}
\]
where in the last step we used~\eqref{Gmn}. This together with~\eqref{boundedgrowth1}, \eqref{nm-1} and~\eqref{6:10} yields that 
\[
\begin{split}
\|h_k(x)-x\|_k &\le DK\left (\frac{\mu_k}{\mu_{\lfloor\widetilde\mu^{-1}(\eta_{n})\rfloor+1}}\right )^a+ KM\sum_{j=\lfloor\widetilde\mu^{-1}(\eta_{n})\rfloor+1}^{k-1}\left (\frac{\mu_k}{\mu_{j+1}}\right )^a \frac{\mu_j'}{\mu_j} \\
&\le DK\left (\frac{\mu_{\lfloor\widetilde\mu^{-1}(\eta_{n+1})\rfloor}}{\mu_{\lfloor\widetilde\mu^{-1}(\eta_{n})\rfloor+1}}\right )^a +KM\left (\frac{\mu_{\lfloor\widetilde\mu^{-1}(\eta_{n+1})\rfloor}}{\mu_{\lfloor\widetilde\mu^{-1}(\eta_{n})\rfloor+1}}\right )^a\sum_{j=\lfloor\widetilde\mu^{-1}(\eta_{n})\rfloor+1}^{k-1}\frac{\mu_j'}{\mu_j}\\
&\le DKe^a+\theta KMe^a.
\end{split}
\]
Hence, the first estimate in~\eqref{Id+} holds. Similarly, one can establish the second. The proof of the theorem is completed.
\end{proof}

In order to give an application of Theorem~\ref{NTL} to nonuniform behavior, we recall the concept of strong nonuniform $\mu$-dichotomy.
\begin{definition}\label{s-n-d}
    We say that~\eqref{eq:Lin} admits a strong nonuniform $\mu$-dichotomy if the operators $A_n$ are invertible and there exists a sequence of projections $\{P_k\}_{k \in \mathbb{N}}$ and constants $N\ge 1$, $\tilde{\nu}\ge \nu>0$ and $\varepsilon \ge 0$ such that (od\ref{def:od1}) in Definition~\ref{def:ordinary-dich} is satisfied and in addition,
    \[
    \|\Phi_A(m, k)P_k\| \le N\left (\frac{\mu_m}{\mu_k}\right )^{-\nu} \mu_k^\varepsilon \quad \text{and} \quad \|\Phi_A(m, k)(\Id-P_k)\| \le N \left (\frac{\mu_m}{\mu_k}\right )^{\tilde \nu}\mu_k^\varepsilon,
    \]
    for $m\ge k$, and 
    \[
    \|\Phi_A(m, k)(\Id-P_k)\| \le N\left (\frac{\mu_k}{\mu_m}\right )^{-\nu} \mu_k^\varepsilon \quad \text{and} \quad \|\Phi_A(m, k)P_k\| \le N\left (\frac{\mu_k}{\mu_m}\right )^{\tilde \nu}\mu_k^\varepsilon,
    \]
    for $m\le k$.
\end{definition}
\begin{corollary}\label{cor-NTL}
    Assume that the following conditions hold:
\begin{enumerate}
\item \eqref{eq:Lin} admits a strong nonuniform $\mu$-dichotomy. 
\item there exists $M>0$ such that 
\begin{equation}\label{BND-new}
    \|g_n(x)\| \le M \mu_n' \mu_n^{-(1+\varepsilon)}\quad \text{for $n\in \mathbb N$ and $x\in X$,}
\end{equation}
where $\mu_n'=\mu_{n+1}-\mu_n$ and with $\varepsilon\ge 0$ as in Definition~\ref{s-n-d};
\item  there exists $c>0$ such that 
\begin{equation}\label{LIP-new}
    \|g_n(x)-g_n(y)\|\le c\mu_n' \mu_n^{-(1+\varepsilon)}\|x-y\|, \quad \text{for $n\in \N$ and $x,y\in X$.}
\end{equation}
\end{enumerate}
Then, provided that $c$ is sufficiently small, there exists a sequence $(h_n)_{n\in \N}$ of homeomorphisms $h_n\colon X \to X$ such that the following holds:
\begin{enumerate}
\item $h_{n+1}\circ (A_n+g_n)=A_n \circ h_n$ for $n\in \N$;
\item there exists $D>0$ such that for $x\in X$ and $n\in \N$,
\[
\|h_n(x)-x\| \le D \quad \text{and} \quad \|h_n^{-1}(x)-x\| \le D;
\]
\item there exist $\bar K, \varrho>0$ and $\delta>0$ such that 
\[
\|h_n(x)-h_n(y)\| \le \bar K\mu_n^{\varrho\varepsilon}\|x-y\|^\varrho \quad \text{and} \quad \|h_n^{-1}(x)-h_n^{-1}(y)\| \le \bar K\mu_n^{\varrho \varepsilon}\|x-y\|^{\varrho},
\]
for $n\in \N$ and $x, y\in X$ with $\|x-y\|\le \delta \mu_n^{-\varepsilon}$.
\end{enumerate}
\end{corollary}

\begin{proof}
For $n\in \N$ and $x\in X$, let 
\[
\|x\|_n:=\|x\|_n^s+\|x\|_n^u,
\]
where
\[
\|x\|_n^s:=\sup_{m\ge n} \left (\|\Phi_A(m, n)P_nx\| \left (\frac{\mu_m}{\mu_n}\right )^\nu \right )+\sup_{m< n}\left (\|\Phi_A (m, n)P_nx\| \left (\frac{\mu_n}{\mu_m}\right )^{-\tilde \nu}\right )
\]
and 
\[
\begin{split}
\|x\|_n^u &:=\sup_{m\le n}\left (\|\Phi_A(m, n)(\Id-P_n)x\| \left (\frac{\mu_n}{\mu_m}\right )^{\nu}\right )\\
&\phantom{=}+\sup_{m>n}
\left (\|\Phi_A(m, n)(\Id-P_n)x\| \left (\frac{\mu_m}{\mu_n}\right )^{-\tilde \nu}\right ).
\end{split}
\]
Clearly, 
\begin{equation}\label{lyapnorms}
\|x\| \le \|x\|_n \le 4N\mu_n^\varepsilon \|x\|, \quad \text{for $x\in X$ and $n\in \N$.} 
\end{equation}
By arguing as in the proof of~\cite[Theorem 5.2]{SilvaCPAA} one can show that~\eqref{eq:Lin} admits $\mu$-dichotomy with respect to the sequence of norms $\left \{\| \cdot \|_n\right \}_{n\in \N}$, and that~\eqref{boundedgrowth1} and~\eqref{boundedgrowth2} hold with some $K>0$ and $a=\tilde \nu>0$.

Observe that it follows from~\eqref{eq:bound-rates-1-a}, \eqref{BND-new} and~\eqref{lyapnorms} that 
\[
\|g_n(x)\|_{n+1}\le 4N\mu_{n+1}^\varepsilon \|g_n(x)\|\le  4N\theta^\varepsilon  \mu_n^\varepsilon \|g_n(x)\| \le 4NM\theta^\varepsilon \frac{\mu_n'}{\mu_n}, 
\]
for $x\in X$ and $n\in \N$.
Hence, \eqref{BND} holds. Similarly,  using~\eqref{eq:bound-rates-1-a}, \eqref{LIP-new} and~\eqref{lyapnorms} we have
\[
\|g_n(x)-g_n(y)\|_{n+1} \le 4Nc\theta^\varepsilon \frac{\mu_n'}{\mu_n} \|x-y\|_n \quad \text{for $x, y\in X$ and $n\in \N$,}
\]
yielding~\eqref{LIP}. The conclusions of the corollary now follow readily from Theorem~\ref{NTL}. For example, \eqref{Id+} and~\eqref{lyapnorms} give that 
\[
\|h_n(x)-x\| \le \|h_n(x)-x\|_n \le D, \quad \text{for $x\in X$ and $n\in \N$.}
\]
Moreover, if $\|x-y\| \le \frac{\delta}{4N}\mu_n^{-\varepsilon}$, then it follows from~\eqref{lyapnorms} that 
$\|x-y\|_n \le \delta$, and consequently by~\eqref{Holdercont} we have that 
\[
\|h_n(x)-h_n(y)\| \le \|h_n(x)-h_n(y)\|_n \le \tilde K\|x-y\|_n^\varrho \le \tilde K (4N)^\varrho \mu_n^{\varrho \varepsilon} \|x-y\|^\varrho.
\] Similarly, one can establish the desired estimates involving $h_n^{-1}$.
\end{proof}

\begin{remark}\label{comparison}
    We note that the result similar to Corollary~\ref{cor-NTL} was formulated in~\cite[Theorem 3]{BV-09}. However, there are important differences between these two results. Namely, we deal with one-sided dynamics while in~\cite{BV-09} the authors considered the two-sided case. More importantly, in~\cite{BV-09} the authors deal with growth rates of the form $\mu_n=e^{\bar \mu_n}$, with $\bar \mu_n$ being an integer for each $n$ (see the bottom of~\cite[p.1979]{BV-09}). Consequently, the case of polynomial behavior is not covered by the results in~\cite{BV-09}, contrary to what is claimed in the abstract of~\cite{BV-09}.
\end{remark}

\section{Nonautonomous smooth linearization: A Sternberg-type theorem}
We begin by introducing some notation. Let $k\in \N$ and $\mu=(\mu_n)_{n\in \N_0}$ be a growth rate. By $\mathcal O_\mu^k$ we denote the set of all sequences $(f_n)_{n\in \N}$ of maps $f_n\colon \mathbb R^d \to \mathbb R^d$ of class $C^k$ such that:
\begin{itemize}
    \item for $n\in \N$, $f_n(0)=0$ and $Df_n(0)=0$;
    \item there exists $M>0$ such that 
    \begin{equation}\label{ok}
        \|D^j f_n(x)\| \le M\frac{\mu_n'}{\mu_n}, \quad \text{for $n\in \N$, $x\in \mathbb R^d$ and $0\le j\le k$.}
    \end{equation}
\end{itemize}
In the particular case, when $\mu_n=e^n$, $n\in \N_0$ we will write $\mathcal O^k$ instead of $\mathcal O_\mu^k$. Furthermore, $\langle \cdot, \cdot \rangle$ will denote the standard Euclidean scalar product, and $B_r(0)$ is an open ball in $\R^d$ of radius r centered in $0$.

Our goal is to prove the following result.
\begin{theorem}\label{Sieg}
Let $\mu=(\mu_n)_{n\in \N_0}$ be a growth rate that satisfies~\eqref{eq:bound-rates-1-a} for some $\theta \ge 1$ and 
$\mathbb A=\{A_n\}_{n\in \mathbb N}$  a sequence of invertible linear operators on $\mathbb R^d$ such that~\eqref{eq:Lin} admits a $\mu$-dichotomy.  Furthermore, suppose that there exist $K, a>0$ satisfying~\eqref{bg-1} and~\eqref{bg-2} and  that 
\[
\Sigma_{\mu D, \mathbb A}=[a_1, b_1]\cup \ldots \cup [a_r, b_r],
\]
with $1\le r\le d$ and 
\begin{equation}\label{abi}
a_1\le b_1 <\ldots <a_r \le b_r.
\end{equation}
Then, for every $\ell \in \N$ there exists $t\in \N$, $t\ge \ell$ with the property that if
\begin{equation}\label{resonance}
[a_i, b_i]\cap [\langle a, q\rangle, \langle b, q\rangle]=\emptyset \quad \text{for $1\le i\le r$ and $q\in \N_0^r$ with $2\le |q|\le t$,}
\end{equation}
where $a=(a_1, \ldots, a_r)$, $b=(b_1, \ldots, b_r)$, $|q|=q_1+\ldots +q_r$, $q=(q_1, \ldots, q_r)$
holds and $(g_n)_{n\in \N}\in \mathcal O_\mu^{t+2}$ with $A_n+g_n$ being  a homeomorphism for each $n$, then there are $p, \tilde p>0$, $r\in (0, p)$, $\tilde r\in (0, \tilde p)$, and for each $k\in \N$ two $C^\ell$ diffeomorphisms $h_k\colon B_r(0)\to h_k(B_r(0))\subset B_{\tilde p}(0)$ and $\bar h_k\colon B_{\tilde r}(0)\to \bar h_k(B_{\tilde r}(0))\subset B_p(0)$ satisfying the following properties:
\begin{enumerate}
\item $\bar h_k(h_k(x))=x$ for each $k\in \N$ and $x\in B_r(0)$ such that $h_k(x)\in B_{\tilde r}(0)$;
\item $h_k(\bar h_k(x))=x$ for each $k\in \N$ and $x\in B_{\tilde r}(0)$ such that $\bar h_k(x)\in B_r(0)$;
\item for $k\in \N$ and $x\in B_r(0)$ such that $A_kx+g_k(x)\in B_{\tilde r}(0)$, 
\[
\bar h_{k+1}(A_k+g_k(x))=A_k h_k(x);
\]
\item 
\[
\lim_{x\to 0}h_k(x)=0 \quad \text{and} \quad \lim_{x\to 0}\bar h_k(x)=0 \quad \text{uniformly in $k$.}
\]
\end{enumerate}
\end{theorem}
\begin{proof}
We consider the sequence $\mathbb B=(B_n)_{n\in \N}$ given by 
\begin{equation}\label{B_n}
B_n:=\Phi_{A}(\lfloor \widetilde{\mu}^{-1}(e^n)\rfloor+1,\lfloor \widetilde{\mu}^{-1}(e^{n-1})\rfloor+1), \quad n\in \N.
\end{equation}
 By Corollary~\ref{cor-442} we have that $\mathbb B$ admits an exponential dichotomy and that~\eqref{supp} holds. Moreover, Theorem~\ref{dicspec} gives that $\Sigma_{ED, \mathbb B}=\Sigma_{\mu D, \mathbb A}$. For $n\in \N$, let  $f_n\colon X\to X$ be given by~\eqref{fnx} with $\mathcal G(m, n)$ as in~\eqref{gmn}.
We now claim that $(f_n)_n\in \mathcal O^{t+2}$. Clearly, each $f_n$ is of class $C^{t+2}$. Moreover, since $g_n(0)=0$ for each $n\in \N$, we have that $\mathcal G(m, n)(0)=0$ for $m, n\in \N$. This gives that $f_n(0)=0$ for each $n\in \N$. In addition, 
\[
Df_n(0)=\sum_{j=\lfloor \widetilde{\mu}^{-1}(\eta_{n-1})\rfloor+1}^{\lfloor \widetilde{\mu}^{-1}(\eta_{n})\rfloor}\Phi_A(\lfloor \widetilde{\mu}^{-1}(\eta_{n})\rfloor+1, j+1) Dg_j (0)D\mathcal G(j, \lfloor \widetilde{\mu}^{-1}(\eta_{n-1})\rfloor+1)(0)=0
\]
for $n\in \N$,
since $Dg_n(0)=0$ for $n\in \N$. Proceeding as in the proof of~\eqref{BND-3} we find that there exists $M_0>0$ such that 
\[
\|f_n(x)\| \le M_0, \quad \text{for $n\in \N$ and $x\in \R^d$.}
\]
Next, we claim  that there exists $\tilde a \ge a$ such that 
\begin{equation}\label{krt}
\|D\mathcal G(m, n)(x)\| \le  K \left (\frac{\mu_m}{\mu_n}\right )^{\tilde a}, \quad \text{for $m\ge n$ and $x\in \mathbb R^d$.}
\end{equation}
To this end, we observe (see~\eqref{gj}) that for $m\ge n$ and $x\in \R^d$,
\[
D\mathcal G(m, n)(x)=\Phi_A(m, n)+\sum_{j=n}^{m-1}\Phi_A(m, j+1)Dg_j(\mathcal G(j, n)(x))D\mathcal G(j, n)(x),
\]
and thus
\[
\begin{split}
\|D\mathcal G(m, n)(x)\| &\le K\left (\frac{\mu_m}{\mu_n}\right )^a+KM\sum_{j=n}^{m-1}\left (\frac{\mu_m}{\mu_{j+1}}\right )^a \frac{\mu_j'}{\mu_j}\|D\mathcal G(j, n)(x)\| \\
&\le K\left (\frac{\mu_m}{\mu_n}\right )^a+KM\sum_{j=n}^{m-1}\left (\frac{\mu_m}{\mu_{j}}\right )^a \frac{\mu_j'}{\mu_j}\|D\mathcal G(j, n)(x)\|.
\end{split}
\]
Hence, 
\[
\left (\frac{\mu_n}{\mu_m}\right )^a \|D\mathcal G(m, n)(x)\| \le K+KM \sum_{j=n}^{m-1}\frac{\mu_j'}{\mu_j}\left (\frac{\mu_n}{\mu_j}\right )^a \|D\mathcal G(j, n)(x)\|,
\]
which together with Gronwall's lemma and~\eqref{nm-1} yields  
\[
\|D\mathcal G(m, n)(x)\| \le K\left (\frac{\mu_m}{\mu_n}\right )^a e^{KM\sum_{j=n}^{m-1}\frac{\mu_j'}{\mu_j}} \le K\left (\frac{\mu_m}{\mu_n}\right )^{a+KM\theta}.
\]
Therefore, \eqref{krt} holds with $\tilde a:=a+KM\theta$.

Since 
\[
\begin{split}
&Df_n(x)=\\
&\sum_{j=\lfloor \widetilde{\mu}^{-1}(\eta_{n-1})\rfloor+1}^{\lfloor \widetilde{\mu}^{-1}(\eta_{n})\rfloor}\Phi_A(\lfloor \widetilde{\mu}^{-1}(\eta_{n})\rfloor+1, j+1)Dg_j (\mathcal G(j, \lfloor \widetilde{\mu}^{-1}(\eta_{n-1})\rfloor+1)(x))D\mathcal G(j, \lfloor \widetilde{\mu}^{-1}(\eta_{n-1})\rfloor+1)(x),
\end{split}
\]
we have using~\eqref{eq:bound2}, \eqref{bg-1}, \eqref{nm-1}, \eqref{ok} and~\eqref{krt}  that 
\begin{equation}\label{dfnx}
\begin{split}
\|Df_n(x)\| &\le K^2 M\sum_{j=\lfloor \widetilde{\mu}^{-1}(\eta_{n-1})\rfloor+1}^{\lfloor \widetilde{\mu}^{-1}(\eta_{n})\rfloor}\left (\frac{\mu_{\lfloor \widetilde{\mu}^{-1}(\eta_{n})\rfloor+1}}{\mu_{j+1}}\right )^a  \frac{\mu_j'}{\mu_j} \left (\frac{\mu_j}{\mu_{\lfloor \widetilde{\mu}^{-1}(\eta_{n-1})\rfloor+1}}\right)^{\tilde a} \\
&\le K^2M\left (\frac{\mu_{\lfloor \widetilde{\mu}^{-1}(\eta_{n})\rfloor+1}}{\mu_{\lfloor \widetilde{\mu}^{-1}(\eta_{n-1})\rfloor+1}}\right )^{a+\tilde a}\sum_{j=\lfloor \widetilde{\mu}^{-1}(\eta_{n-1})\rfloor+1}^{\lfloor \widetilde{\mu}^{-1}(\eta_{n})\rfloor}\frac{\mu_j'}{\mu_j} \\
&\le \theta K^2M\left (\frac{\mu_{\lfloor \widetilde{\mu}^{-1}(\eta_{n})\rfloor+1}}{\mu_{\lfloor \widetilde{\mu}^{-1}(\eta_{n-1})\rfloor+1}}\right )^{a+\tilde a}\log \left (\frac{\mu_{\lfloor \widetilde{\mu}^{-1}(\eta_{n})\rfloor+1}}{\mu_{\lfloor \widetilde{\mu}^{-1}(\eta_{n-1})\rfloor+1}}\right ) \\
&\le \theta^{1+2(a+\tilde a)}K^2M e^{a+\tilde a}\log (\theta^2 e).
\end{split}
\end{equation}
Hence, we find that there exists $M_1>0$ such that 
\[
\|Df_n(x)\| \le M_1, \quad \text{for $n\in \N$ and $x\in \R^d$.}
\]
Proceeding in the same manner, we find that for each $0\le j \le t+2$, there is $M_j>0$ such that 
\[
\|D^j f_n(x)\| \le M_j, \quad \text{for $n\in \N$ and $x\in \mathbb R^d$.}
\]
This implies that $(f_n)_{n\in \N}\in \mathcal O^{t+2}$. We are now in a position to apply Corollary~\ref{siegm-cor} to $\mathbb B$ and $(f_n)_{n\in \N}$. Hence, there are $p, \tilde p>0$, $r\in (0, p)$, $\tilde r\in (0, \tilde p)$, and for each $k\in \N$ two $C^\ell$ diffeomorphisms $\psi_k\colon B_r(0)\to \psi_k(B_r(0))\subset B_{\tilde p}(0)$ and $\bar \psi_k\colon B_{\tilde r}(0)\to \bar \psi_k(B_{\tilde r}(0))\subset B_p(0)$ with the following properties:
\begin{enumerate}
\item $\bar \psi_k(\psi_k(x))=x$ for each $k\in \N$ and $x\in B_r(0)$ such that $\psi_k(x)\in B_{\tilde r}(0)$;
\item $\psi_k(\bar \psi_k(x))=x$ for each $k\in \N$ and $x\in B_{\tilde r}(0)$ such that $\bar \psi_k(x)\in B_r(0)$;
\item for $k\in \N$ and $x\in B_r(0)$ such that $B_kx+f_k(x)\in B_{\tilde r}(0)$, 
\[
\bar \psi_{k+1}(B_k+f_k(x))=B_k \psi_k(x);
\]
\item 
\[
\lim_{x\to 0}\psi_k(x)=0 \quad \text{and} \quad \lim_{x\to 0}\bar \psi_k(x)=0 \quad \text{uniformly in $k$.}
\]
\end{enumerate}
Similarly to~\eqref{krt}, one can show that 
\begin{equation}\label{krt2}
\|D\mathcal G(m, n)(x)\| \le  K \left (\frac{\mu_n}{\mu_m}\right )^{\tilde a}, \quad \text{for $m\le n$ and $x\in \mathbb R^d$.}
\end{equation}
Since $\mathcal G(m, n)(0)=0$, it follows from~\eqref{krt2} that 
\begin{equation}\label{gmnx}
\|\mathcal G(m, n)(x)\| \le K\left (\frac{\mu_n}{\mu_m}\right )^{\tilde a}\|x\|, \quad \text{for $m\le n$ and $x\in \mathbb R^d$.}
\end{equation}
This implies that there exists $r_1>0$ such that 
\[
\|\mathcal G(\lfloor\widetilde\mu^{-1}(\eta_{n})\rfloor+1, k)(x)\| < r,
\]
for $k\in \N$, $x\in \R^d$ with $\|x\|<r_1$ and where $n\in \N_0$ satisfies~\eqref{nn+1}.   Indeed, we can take $r_1=\frac{r}{Ke^{\tilde a}}$. 
For $k\in \N$, we define $h_k$ on $B_{r_1}(0)$ by~\eqref{hk} (where $n\in \N_0$  is such that~\eqref{nn+1} holds). Clearly, $h_k$ is of class $C^\ell$ on $B_{r_1}(0)$, and since~\eqref{bg-2} implies uniform bound for $\|\Phi_A (k, \lfloor\widetilde\mu^{-1}(\eta_{n})\rfloor+1)\|$ , we have that there is $\tilde p_1>0$ such that $h_k(B_{r_1}(0))\subset B_{\tilde p_1}(0)$. 

We now define $\bar h_k$ for $k\in \N$. By~\eqref{bg-2}, there exists $\tilde r_1>0$ such that 
\[
\|\Phi_A(\lfloor \widetilde{\mu}^{-1}(\eta_{n+1})\rfloor+1, k)x\|<\tilde r,
\]
for $k\in \N$, $x\in \R^d$ with $\|x\|<\tilde r_1$, where $n\in \N_0$ satisfies~\eqref{nn+1}. We now define
\[
\bar h_k:= \mathcal G(k, \lfloor\widetilde\mu^{-1}(\eta_{n})\rfloor+1) \circ \bar \psi_{n+1} \circ \Phi_A ( \lfloor\widetilde\mu^{-1}(\eta_{n})\rfloor+1, k),
\]
where $n\in \N_0$ is such that~\eqref{nn+1} holds. 
From the preceding discussion, we see that $\bar h_k$ is well defined and of class $C^\ell$ on $B_{\tilde r_1}(0)$ for each $k\in \N$. From~\eqref{krt} 
we find that there exists $p_1>0$ such that $\bar h_k(B_{\tilde r_1}(0))\subset B_{p_1}(0)$. It is straightforward to verify that the maps $h_k$ and $\bar h_k$ have the desired properties.
\end{proof}

\section{Appendix}
In this section, we establish two results that were used in previous sections.
\subsection{Nonautonomous topological linearization on the half-line}
We recall the following result which is established in the proof of~\cite[Theorem 5]{BV}.
\begin{theorem}\label{5-12-BV}
Let $\{\|\cdot \|_n\}_{n\in \mathbb Z}$ be a sequence of norms on $X$ equivalent to $\| \cdot \|$, $\{A_n\}_{n\in \mathbb Z}$ a sequence of invertible bounded operators on $X$ which admits an exponential dichotomy with respect to the sequence of norms $\{\|\cdot \|_n\}_{n\in \mathbb Z}$, and with the property that there exists $M>0$ such that 
\begin{equation}\label{bg-new}
\|A_n x\|_{n+1}\le M\|x\|_n \quad \text{and} \quad \|A_n^{-1}x\|_n \le M\|x\|_{n+1},
\end{equation}
for $n\in \mathbb Z$ and $x\in X$. Moreover, let $\{f_n\}_{n\in \mathbb Z}$ be a sequence of maps $f_n\colon X\to X$ such that there exist $c, C>0$ such that:
\begin{enumerate}
    \item for $x\in X$ and $n\in \mathbb Z$, $\|f_n(x)\|_{n+1}\le C\|x\|_n$;
    \item for $x, y\in X$ and $n\in \mathbb Z$,
    \[
    \|f_n(x)-f_n(y)\|_{n+1}\le c\|x-y\|_n.
    \]
\end{enumerate}
Then, provided that $c$ is sufficiently small, there exists a sequence $\psi_n\colon X\to X$, $n\in \mathbb Z$ satisfying
\[
\psi_{n+1}\circ (A_n+f_n)=A_n \circ \psi_n, \quad \text{for $n\in \mathbb Z$.}
\]
In addition, there are $D, \varrho>0$ such that 
\[
\|\psi_n(x)-\psi_n(y)\|_n\le D\|x-y\|_{n}^\varrho
\]
and
\[
\|\psi_n^{-1}(x)-\psi_n^{-1}(y)\|_n\le D\|x-y\|_{n}^\varrho,
\]
for $n\in \mathbb Z$ and $x, y\in X$ such that $\|x-y\|_{n}\le 1$. Finally, 
\[
\|\psi_n(x)-x\|_{n}\le D \quad \text{and} \quad \|\psi_n^{-1}(x)-x\|_{n}\le D, 
\]
for $n\in \mathbb Z$ and $x\in X$.
\end{theorem}
We now establish the version of Theorem~\ref{5-12-BV} for one-sided dynamics.
\begin{corollary}\label{BVcor}
    Let $\{\|\cdot \|_n\}_{n\in \N}$ be a sequence of norms on $X$ equivalent to $\| \cdot \|$, $\{A_n\}_{n\in \N}$ a sequence of invertible bounded operators on $X$ that admits an exponential dichotomy with respect to the sequence of norms $\{\|\cdot \|_n\}_{n\in \N}$, and with the property that there exists $M>0$ such that~\eqref{bg-new} holds for $n\in \N$ and $x\in X$.
 Moreover, let $\{f_n\}_{n\in \N}$ be a sequence of maps $f_n\colon X\to X$ such that there exist $c, C>0$ such that:
\begin{enumerate}
    \item for $x\in X$ and $n\in \N$, $\|f_n(x)\|_{n+1}\le C\|x\|_n$;
    \item for $x, y\in X$ and $n\in \N$,
    \[
    \|f_n(x)-f_n(y)\|_{n+1}\le c\|x-y\|_n.
    \]
\end{enumerate}
Then, provided that $c$ is sufficiently small, there exists a sequence $\psi_n\colon X\to X$, $n\in \N$ satisfying
\[
\psi_{n+1}\circ (A_n+f_n)=A_n \circ \psi_n, \quad \text{for $n\in \N$.}
\]
In addition, there are $D, \varrho>0$ such that 
\[
\|\psi_n(x)-\psi_n(y)\|_n\le D\|x-y\|_{n}^\varrho
\]
and
\[
\|\psi_n^{-1}(x)-\psi_n^{-1}(y)\|_n\le D\|x-y\|_{n}^\varrho,
\]
for $n\in \mathbb N$ and $x, y\in X$ such that $\|x-y\|_{n}\le 1$. Finally, 
\[
\|\psi_n(x)-x\|_{n}\le D \quad \text{and} \quad \|\psi_n^{-1}(x)-x\|_{n}\le D, 
\]
for $n\in \N$ and $x\in X$.
\end{corollary}
\begin{proof}
Let $\|\cdot \|_n=\|\cdot \|$ for $n\le 0$. Moreover, we choose an arbitrary invertible operator $A\in \mathcal B(X)$ which is hyperbolic (i.e. its spectrum does not intersect the unit circle), and such that 
\[
\text{Im} P_1=\{x\in X: \ \lim_{n\to \infty}A^n x=0\} \quad \text{and} \quad  \Ker P_1=\{x\in X: \ \lim_{n\to \infty}A^{-n}x=0\},
\]
where $\{P_n\}_{n\in \N}$ is a sequence of projections associated with the dichotomy of $\{A_n\}_{n\in \N}$.
Set $A_n:=A$ for $n\le 0$.  It is straightforward to verify that the sequence $\{A_n\}_{n\in \mathbb Z}$ admits an exponential dichotomy with respect to the sequence of norms $\{\|\cdot \|_n\}_{n\in \mathbb Z}$. Moreover, \eqref{bg-new} holds. Finally, we set $f_n=0$ for $n\le 0$. Obviously, the sequence $\{f_n\}_{n\in \mathbb Z}$ satisfies the assumptions of Theorem~\ref{5-12-BV}. The conclusion of the corollary now follows readily from Theorem~\ref{5-12-BV}.
\end{proof}

\subsection*{Nonautonomous smooth linearization on half-line}
We recall the following result, which is essentially established in~\cite{CDS} (see Remark~\ref{comments}).
\begin{theorem}\label{sternberg}
Let $\mathbb B=\{B_n\}_{n\in \mathbb Z}$ be a sequence of invertible linear operators on $\mathbb R^d$ that admits an  exponential dichotomy and with the property that 
\begin{equation}\label{supp}
    \sup_{n\in \mathbb Z}\|B_n\|<+\infty \quad \text{and}\quad \sup_{n\in \mathbb Z}\|B_n^{-1}\|<+\infty. 
\end{equation}
Suppose that
\[
\Sigma_{ED, \mathbb B}=[a_1, b_1]\cup \ldots \cup [a_r, b_r],
\]
with $1\le r\le d$ and $a_i, b_i$ as in~\eqref{abi}.
Then, for every $\ell \in \N$ there exists $t\in \N$, $t\ge \ell$ with the following property: if~\eqref{resonance} holds  and $(f_n)_{n\in \mathbb Z}\in \mathcal O^{t+2}$ with $B_n+f_n$ being 
a homeomorphism on $\R^d$ for each $n$, then there are $p, \tilde p>0$, $r\in (0, p)$, $\tilde r\in (0, \tilde p)$, and for each $n\in \mathbb Z$ two $C^\ell$ diffeomorphisms $h_n\colon B_r(0)\to h_n(B_r(0))\subset B_{\tilde p}(0)$ and $\bar h_n\colon B_{\tilde r}(0)\to \bar h_n(B_{\tilde r}(0))\subset B_p(0)$ with the following properties:
\begin{enumerate}
\item $\bar h_n(h_n(x))=x$ for each $n\in \mathbb Z$ and $x\in B_r(0)$ such that $h_n(x)\in B_{\tilde r}(0)$;
\item $h_n(\bar h_n(x))=x$ for each $n\in \mathbb Z$ and $x\in B_{\tilde r}(0)$ such that $\bar h_n(x)\in B_r(0)$;
\item for $n\in \mathbb Z$ and $x\in B_r(0)$ such that $B_nx+f_n(x)\in B_{\tilde r}(0)$, 
\[
\bar h_{n+1}(B_nx+f_n(x))=B_n h_n(x);
\]
\item 
\[
\lim_{x\to 0}h_n(x)=0 \quad \text{and} \quad \lim_{x\to 0}\bar h_n(x)=0 \quad \text{uniformly in $n$.}
\]
\end{enumerate}
\end{theorem}
\begin{remark}\label{comments}
We note that in~\cite{CDS} the authors have formulated a continuous-time version of Theorem~\ref{Sieg}.
However, the proof proceeds by discretization of time, thus essentially yielding the version for discrete time as a by-product.  Moreover, $(f_n)_{n\in \mathbb Z}\in \mathcal O^{t+2}$ means that $f_n(0)=0$ and $Df_n(0)=0$ for $n\in \mathbb Z$, and there exists $M>0$ such that 
\[
\|D^j f_n(x)\| \le M, \quad \text{for $n\in \mathbb Z$, $x\in \mathbb R^d$ and $0\le j\le t+2$.}
\]
\end{remark}
We now establish the version of Theorem~\ref{Sieg} for a half-line.
\begin{corollary}\label{siegm-cor}
    Let $\mathbb B=\{B_n\}_{n\in \N}$ be a sequence of invertible linear operators on $\mathbb R^d$ that admits an  exponential dichotomy and with the property that 
\[
    \sup_{n\in \N}\|B_n\|<+\infty \quad \text{and}\quad \sup_{n\in \N}\|B_n^{-1}\|<+\infty. 
\] 
Suppose that
\[
\Sigma_{ED, \mathbb B}=[a_1, b_1]\cup \ldots \cup [a_r, b_r],
\]
with $1\le r\le d$ and $a_i, b_i$ as in~\eqref{abi}.
Then, for every $\ell \in \N$ there exists $t\in \N$, $t\ge \ell$ with the following property: if~\eqref{resonance} holds and $(f_n)_{n\in \mathbb Z}\in \mathcal O^{t+2}$ with $B_n+f_n$ being 
a homeomorphism on $\R^d$ for each $n$, then there are $p, \tilde p>0$, $r\in (0, p)$, $\tilde r\in (0, \tilde p)$, and for each $n\in \N_0$ two $C^\ell$ diffeomorphisms $h_n\colon B_r(0)\to h_n(B_r(0))\subset B_{\tilde p}(0)$ and $\bar h_n\colon B_{\tilde r}(0)\to \bar h_n(B_{\tilde r}(0))\subset B_p(0)$ with the following properties:
\begin{enumerate}
\item $\bar h_n(h_n(x))=x$ for each $n\in \N$ and $x\in B_r(0)$ such that $h_n(x)\in B_{\tilde r}(0)$;
\item $h_n(\bar h_n(x))=x$ for each $n\in \N_0$ and $x\in B_{\tilde r}(0)$ such that $\bar h_n(x)\in B_r(0)$;
\item for $n\in \N_0$ and $x\in B_r(0)$ such that $B_nx+f_n(x)\in B_{\tilde r}(0)$, 
\[
\bar h_{n+1}(B_nx+f_n(x))=B_n h_n(x);
\]
\item 
\[
\lim_{x\to 0}h_n(x)=0 \quad \text{and} \quad \lim_{x\to 0}\bar h_n(x)=0 \quad \text{uniformly in $n$.}
\]
\end{enumerate}
\end{corollary}

\begin{proof}
As shown in the proof of~\cite[Theorem 1]{BDZ}, we can extend $\mathbb B$ to a two-sided sequence $\mathbb B'=(B_n')_{n\in \mathbb Z}$ of invertible operators such that $B_n=B_n'$ for $n\in \N$ and
$\Sigma_{ED, \mathbb B}=\Sigma_{ED,\mathbb B'}$ (in particular, $\mathbb B'$ admits an exponential dichotomy). Moreover, 
\[
\sup_{n\in \mathbb Z}\|B_n'\|<+\infty \quad \text{and} \quad \sup_{n\in \mathbb Z}\|(B_n')^{-1}\|<+\infty.\]
Set $f_n:=0$ for $n\le 0$. Applying Theorem~\ref{sternberg} to $\mathbb B'$ and $(f_n)_{n\in \mathbb Z}$ immediately yields the desired conclusion.
\end{proof}

\end{document}